\documentclass[]{article}
\usepackage[]{geometry}
\usepackage{graphicx,subfigure}
\usepackage{amsmath, amsfonts, amsthm, amssymb, amstext, multicol, wrapfig, caption, , float, multirow, tabularx,bm,ragged2e,array,scrextend, verbatim,enumerate,mathtools,indentfirst}
\usepackage[colorlinks]{hyperref}
\hypersetup{
	citecolor = blue
}
\usepackage{authblk}
\usepackage{xcolor}

\newcommand{\id}{\mathrm d}
\newcommand{\vc}{\bm}
\newcommand{\bs}{\boldsymbol}

\newcommand{\gr}{\mbox{g}}

\renewcommand{\tilde}{\widetilde}
\newcommand{\pard}[2]{\frac{\partial #1}{\partial #2}}

\newcommand{\R}{\mathbb{R}}

\newtheorem{rem}{Remark}

\title{Enforcing conserved quantities in Galerkin truncation and finite volume discretization}
\author[]{Zachary T. Hilliard} 
\author[]{Mohammad Farazmand\thanks{Corresponding author: farazmand@ncsu.edu}}
\affil[]{Department of Mathematics, North Carolina State University, 2311 Stinson Drive, Raleigh, NC 27695, USA}

\date{}

\begin{document}

\maketitle

\begin{abstract}
Finite-dimensional truncations are routinely used to approximate partial differential equations (PDEs), either to obtain numerical solutions or to derive reduced-order models. The resulting discretized equations are known to violate certain physical properties of the system. In particular, first integrals of the PDE may not remain invariant after discretization. Here, we use the method of reduced-order nonlinear solutions (RONS) to ensure that the conserved quantities of the PDE survive its finite-dimensional truncation. In particular, we develop two methods: Galerkin RONS and finite volume RONS. Galerkin RONS ensures the conservation of first integrals in Galerkin-type truncations, whether used for direct numerical simulations or reduced-order modeling. Similarly, finite volume RONS conserves any number of first integrals of the system, including its total energy, after finite volume discretization. Both methods are applicable to general time-dependent PDEs and can be easily incorporated in existing Galerkin-type or finite volume code. We demonstrate the efficacy of our methods on two examples: direct numerical simulations of the shallow water equation and a reduced-order model of the nonlinear Schr\"odinger equation. As a byproduct, we also generalize RONS to phenomena described by a system of PDEs.
\end{abstract}

\section{Introduction}\label{sec:intro}
Galerkin and Petrov-Galerkin methods are routinely used to obtain numerical solutions to partial differential equations (PDEs) and to obtain reduced-order models of the full PDE~\cite{karniadakis2013,Rowley2017}.
Galerkin-type methods approximate the true solution $u(\vc x,t)$ of the PDE with a finite-dimensional truncation of the form,
\begin{equation}
\hat{u}(\vc x,t) = \sum_{i=1}^N a_i(t)\phi_i(\vc x),
\end{equation}
where  $\phi_i$ are some prescribed modes and $a_i$ are the corresponding amplitudes.
The approximate solution $\hat u$ evolves on the $N$-dimensional linear subspace spanned by the modes $\{\phi_i\}_{i=1}^N$. Galerkin-type truncations lead to a set of ordinary differential equations (ODEs) for the amplitudes $a_i(t)$, whose solution dictates the evolution of the approximate solution $\hat u$.

There are two main regimes in which Galerkin truncations are used: direct numerical simulations (DNS) and reduced-order modeling. In DNS, the goal is to approximate the solution of the PDE with high accuracy and therefore the dimension of the linear subspace $N$ is chosen to be as large as possible. Depending on the PDE, the choice of the modes $\phi_i$ can include Fourier modes, finite elements, piece-wise polynomials, etc. 

Many applications, such as inverse problems, sensitivity analysis, control, and uncertainty quantification, require repeated solves of the PDE. In this case, using the costly DNS solution may be impractical; therefore reduced-order models are used instead.
For reduced-order modeling purposes, one seeks an approximate solution with as few modes $N$ as possible, to reduce the computational cost. In other words, one sacrifices accuracy for computational efficiency. Nonetheless, if the modes $\phi_i$ are chosen carefully and the corresponding amplitudes $a_i$ are evolved appropriately, the reduced-order solution will have sufficient accuracy to reproduce many relevant features of the true solution. Several equation-driven and data-driven methods, such as proper orthogonal decomposition (POD)~\cite{berkooz93}, dynamic mode decomposition (DMD)~\cite{schmid10}, and spectral submanifolds (SSM)~\cite{Haller2016,Szalai2017}, have been developed to discover the appropriate set of modes $\phi_i$ for a particular PDE.

Whether used for direct numerical simulations or for reduced-order modeling, Galerkin-type methods are known to violate certain physical properties of the system~\cite{Bowman2013,Carlberg2018,Majda2012,Modin2020}. In particular, quantities which might be conserved under the dynamics of the PDE, often fail to remain conserved under its Galerkin truncation. As we review in section~\ref{sec:rel_work}, several techniques have been developed to address this shortcoming. However, these methods are either developed for a specific PDE or are only applicable to a particular class of PDEs. Here, using the framework of reduced-order nonlinear solutions (RONS)~\cite{Anderson2022a}, we develop a universal method for ensuring that the conserved quantities of the PDE survive its Galerkin truncation. The resulting method is applicable to a wide range of time-dependent PDEs and can be implemented with simple changes to existing code.
Furthermore, we develop two extensions to the theory of RONS, so that our original contributions can be summarized as follows:
\begin{enumerate}
	\item Galerkin RONS: This method ensures that the conserved quantities of the PDE remain conserved after Galerkin truncation.
	\item RONS for systems of PDEs: RONS was originally developed for scalar-valued PDEs. Here, we formulate it for a system of PDEs. In particular, this allows the use of RONS on PDEs with higher-order time derivatives.
	\item Finite volume RONS: Finite volume methods are designed to respect the conservation of the state variable. However, certain conservation-law PDEs admit additional first integrals which may not be respected by the finite volume discretization. To address this issue, we develop finite volume RONS which conserves all known first integrals of the PDE.
\end{enumerate}

\subsection{Related Work}\label{sec:rel_work}
As mentioned earlier, the fact that Galerkin-type truncations do not necessarily respect the physical properties of the PDE is well-known. A number of attempts have been made to address this issue. Here, we review the most relevant studies.

\emph{Symmetry reduction:} By Noeother's theorem, any continuous symmetry of a system leads to a corresponding conserved quantity~\cite{KosmannSchwarzbach2010}. As such, it has long been recognized that symmetries of the system must be preserved in a reduced model~\cite{Marsden1974}. One approach is symmetry reduction whereby, before any further model order reduction, the symmetries of the system are reduced. For instance, the method of slices~\cite{Rowley2000} was developed as a simple computational method for removing translational and rotational symmetries of a system. This method was later generalized to multiple slices to remove the spurious singularities induced by symmetry reduction~\cite{Froehlich2012,Budanur2015,Budanur2017b}. In addition to the method of slices, ideas from differential geometry have been used to reduce general Lie symmetries of a system~\cite{Haller1998b,Marsden1974}.
In any case, removing the symmetries of the system greatly simplifies its dynamics so that subsequent POD-Galerkin model reduction on the symmetry-reduced system succeeds, where it might have failed prior to symmetry reduction~\cite{Mowlavi2018}

\emph{Hamiltonian systems:} Hamiltonian PDEs have a rich geometric structure which may be lost after numerical discretization. Symplectic integrators were developed to ensure that the symplectic two-form associated with a Hamiltonian system is preserved after discretization and consequently the discrete system conserves volume~\cite{vogelaere1956,Ruth1983,Kang1991}. Although first derived for Hamiltonian ODEs, these ideas have been generalized to PDEs~\cite{Bridges2006}. For instance, Modin and Viviani~\cite{Modin2020} developed a structure-preserving discretization of the two-dimensional Euler equation for ideal, incompressible fluids which respects the conservation of its Casimir invariants. Building on these ideas, Cifani et al.~\cite{Cifani2022} showed that 
such structure-preserving discretizations reproduce an energy spectrum which agrees with Kraichnan's theoretical prediction~\cite{kraichnan1967,kraichnan1971}, whereas this prediction is not reproducible by other numerical methods at moderate Reynolds numbers~\cite{Boffetta2010,faraz_cont}.
In addition to DNS, reduced-order modeling methods have been developed that retain the Hamiltonian structure of the system. For example, Mohseni and Peng ~\cite{Peng2016} developed the proper symplectic decomposition (PSD) which ensures that a POD-Galerkin model reduction preserves the symplectic structure of the system. Another example is the symplectic manifold Galerkin (SMG) of Buchfink et. al.~\cite{Buchfink2023} who reduce the model to a lower dimensional manifold as opposed to a linear subspace. 

\emph{Finite volume method:} Many PDEs, which model physical phenomena, are derived from conservation laws. Finite volume methods ensure that these conservation laws are respected by the discretization. By construction, finite volume schemes are conservative in the sense that they conserve the state variables of a system. However, certain conservation-law PDEs, such as the shallow water equation considered in this paper, admit additional first integrals which may not remain conserved after a finite volume discretization. This has motivated the development of \emph{structure-preserving} and \emph{energy}-\emph{preserving} finite volume methods. For example, Fjordholm et. al~\cite{Fjordholm2011} and Tadmor and Zhong~\cite{Tadmor2008} developed structure-preserving schemes which conserve the total energy for the shallow water equation (SWE). Another example appears in Mishra and Tadmor~\cite{Mishra2011} who derived a finite volume scheme which satisfies additional constraints of the system such as divergence-free and curl-free conditions. In the realm of reduced-order modeling, starting from a finite volume discretization, Carlberg et. al \cite{Carlberg2018} developed a new technique to enforce the conservation of the state variables in POD-based Galerkin and Petrov-Galerkin reduced-order models. However, they do not address enforcing any additional first integrals such as the total energy of the system. 

\emph{Stability-preserving methods:} An additional issue with Galerkin-type reduced-order models is their possible loss of stability. More specifically, a stable equilibrium of the system may become unstable after model reduction. Although this issue goes beyond the scope of our work, we review the existing remedies for completeness. Using the control-theoretic idea of balanced truncations from Moore~\cite{Moore1981}, Willcox and Peraire~\cite{Willcox2002} developed balanced POD which restricts the projection onto to the observability and controllability subspaces through empirical approximations of their respective grammians. This method was later enhanced by Rowley~\cite{Rowley2005} who showed that balanced POD can preserve the stability of equilibria. 
In the context of compressible fluid flows, Rowley et al. \cite{Rowley2004} proposed a stability-preserving POD-Galerkin method by using an energy-based inner product instead of the usual $L^2$ inner product. Another example can be found in~\cite{Kalashnikova2014} who studied linear time-invariant systems. They use a full-state feedback controller to force the eigenvalues to lie within the unit disc and ensure the stability of the reduced-order model. These ideas were generalized to linear time-dependent systems in \cite{Mojgani2020} who constrain the singular values of the linear time-varying system through a feedback controller. This procedure limits the transient growth of perturbations and guarantees the stability of the reduced time-dependent system. 

%\emph{Specific methods:} Outside of finite volume schemes, there are also conservative methods for specific PDEs or choice of basis functions. For example, a recent work by Valle and Verstappen \cite{Valle2023} develop an energy preserving scheme for multiphase flows in bubble dynamics. Fei and V\'{a}zquez \cite{Fei1999} proposed two schemes which conserve the total energy for the Sine-Gordon equation. Modin \cite{Modin2020} uses geometric quantization in space and a symplectic integrator in time to preserve the geometric structure of Euler's equation and, consequently, enforce the conservation of the Casimir invariants of the 2D Euler equation for ideal, incompressible fluids. 

The methods reviewed in this section are restrictive in the sense that they either apply to a particular PDE or a particular class of PDEs (e.g., Hamiltonian systems). 
In contrast, our RONS-based method is applicable to any PDE with conserved quantities. Furthermore, the proposed method can be easily implemented on top of existing Galerkin-type and finite volume schemes. Of course, this universality and straightforward implementability come at a cost: although our RONS-based method ensures conservation of the first integrals of the system, it does not necessarily preserve any additional structure such as the symplectic geometry of Hamiltonian systems.
 
\subsection{Outline}
We organize the remainder of the paper in the following manner. In section \ref{sec:math_prelim}, we discuss the general set-up and framework of RONS in the context of Galerkin projections. In section \ref{sec:extensions}, we extend the methodology presented in section \ref{sec:math_prelim} to systems of equations and finite volume discretizations. We present two numerical examples in Section \ref{sec:num_results}. Namely we apply our methodology to the shallow water equation (SWE) in Section \ref{sec:swe} and the nonlinear Schr\"{o}dinger (NLS) equation in section \ref{sec:nls}. Lastly, we present our concluding remarks in Section \ref{sec:conc}.

\section{Mathematical Preliminaries}\label{sec:math_prelim}
RONS is a method for evolving reduced models which depend nonlinearly on a set of time-dependent parameters~\cite{Anderson2022a}. Although it was originally developed for reduced-order modeling, here we use RONS for both model order reduction and DNS. 
In this section, we review the RONS framework in the context of Galerkin-type truncations of scalar, real-valued PDEs of the form,
\begin{equation}\label{eq:gen_pde}
u_t = F(u), \quad u(\vc x,0) = u_0(\vc x),
\end{equation}
where $u:\Omega\times \mathbb{R}^{+} \rightarrow \R $ is the solution of the PDE, $\Omega \subseteq \mathbb{R}^n$ is the spatial domain, $F$ is a differential operator, and $u_0(\vc x)$ is the initial condition of the PDE. In section \ref{sec:RONS_systems}, we extend the theory of RONS to systems of PDEs, but for this review we restrict our attention to scalar equations. 

Additionally, we assume that the PDE has $m$ conserved quantities (or first integrals) denoted by $I_k$, such that 
\begin{equation}\label{eq:consQuantities}
I_k(u(\cdot, t)) = \mbox{const.}, \quad \quad k = 1,2,\ldots,m,
\end{equation}
for all time $t\geq 0$.
The conserved quantities $I_k$ may denote physical quantities such as mass, energy, or momentum, or more abstract quantities that arise from the structure of the PDE such as the Casimir invariants of Euler's equation for ideal fluids \cite{Bowman2013}. 
%An important contribution of RONS is to ensure that reduced-order approximations of the PDE also respect the conserved quantities $I_k$.
We assume that the solution $u$ belongs to an appropriate Hilbert space $H$ endowed with the inner product $\langle \cdot, \cdot\rangle_H$ and the induced norm $\|\cdot\|_H$, and we assume that the boundary conditions for the PDE \eqref{eq:gen_pde} are encoded in the Hilbert space $H$. 

Next, we introduce a Galerkin-type approximate solution $\hat{u}(\vc x,\vc a(t))$ of the form,
\begin{equation}\label{eq:approx_sol}
\hat{u}(\vc x,\vc a(t)) = \sum_{i=1}^{N}a_i(t) \phi_i(\vc x),
\end{equation}
where $N$ is the total number of modes,  $\vc a = (a_1,\cdots, a_N): \mathbb{R}^+\rightarrow \R^N$ is the time-dependent vector of parameters and $\phi_i: \Omega \rightarrow \mathbb{R}$ are the linearly independent basis functions or modes. The basis functions are chosen based on the structure of the PDE \eqref{eq:gen_pde}. For example, the  modes can be Fourier modes or finite elements for high resolution DNS. For reduced-order modeling purposes, $\phi_i$ are a relatively small number of modes, such as those obtained from the proper orthogonal decomposition (POD)~\cite{Sirovich1987,berkooz93} or dynamic mode decomposition (DMD)~\cite{schmid10}. In either case, RONS produces a system of ordinary differential equations (ODEs) to evolve the parameters $\vc a$ in time. 

Even when the solution $u$ is real-valued, the parameters $a_i$ and the modes $\phi_i$ can be complex-valued, as is the case for Fourier modes. For the review in this section, we only consider the case where both the parameters and modes are real-valued. Generalization to complex-valued parameters and modes is straightforward and is included in Appendix \ref{app:rons_complex} for completeness . 

\begin{rem}
We emphasize that the general construction of RONS allows the approximate solution $\hat{u}$ to have a nonlinear dependence on the time-dependent parameters $\vc a(t)$, but we will not consider this feature here. For such nonlinear approximate solutions, we refer the reader to Ref.~\cite{Anderson2022a} which discusses their theoretical underpinnings. Their computational aspects are discussed in Refs.~\cite{Anderson2024a,Anderson2024b}.
\end{rem}

Galerkin truncation can be obtained by substituting the approximate solution~\eqref{eq:approx_sol} into the PDE~\eqref{eq:gen_pde} and truncating the right-hand side to the modes $\{\phi_1,\cdots,\phi_N\}$. An alternative, yet equivalent, approach is to cast Galerkin truncation as an optimization problem~\cite{Anderson2022a}.
More specifically, one seeks to minimize the instantaneous error between the rate of change $\hat{u}_t$ of the approximate solution and the rate $F(\hat{u})$ dictated by the PDE. To this end, we define the cost functional,
\begin{equation}\label{eq:cost_func}
\mathcal{J}(\vc a(t),\dot{\vc a}(t)) = \frac{1}{2} \|\hat{u}_t - F(\hat{u})\|^2_H,
\end{equation}
which quantifies the error between the left- and right-hand sides of the PDE at the approximate solution $\hat u$.
Note that $\hat u_t = \sum_i \dot a_i \phi_i$ depends on $\dot{\vc a}$ and $F(\hat u)$ depends on $\vc a$. As a result, the cost functional $\mathcal J$ depends on both $\vc a(t)$ and $\dot{\vc a}(t)$.

We minimize the function $\mathcal J(\vc a,\dot{\vc a})$ instantaneously in the following sense. At any time $t$, we seek to evolve the parameters $\vc a(t)$ such that $\mathcal J$ is minimized at that instance. More precisely, we solve the instantaneous optimization problem,
\begin{equation}\label{eq:rons_opt}
\min_{\dot{\vc a}\in\mathbb{\R}^N} \mathcal{J}(\vc a,\dot{\vc a}),
\end{equation}
where the minimum is  only taken over $\dot{\vc a}$. Since $\mathcal J$ is quadratic in $\dot{\vc a}$, the optimal solution has an explicit solution given by 
\begin{equation}\label{eq:rons_ode}
\vc M \dot{\vc a} = \vc f(\vc a),
\end{equation}
which is exactly the equation for the standard Galerkin truncation. Here the \textit{metric tensor} $\vc M$, is an $N\times N$ matrix whose elements are given by
\begin{equation}\label{eq:metric_tensor}
M_{ij} =\left\langle\frac{\partial \hat{u}}{\partial a_i},\frac{\partial \hat{u}}{\partial a_j}\right\rangle_H =  \left\langle\phi_i , \phi_j \right\rangle_H, 
\end{equation}
and the vector field $\vc f:\mathbb R^N\to\mathbb R^N$ has components,
\begin{equation}\label{eq:rons_f}
	f_i = \left \langle \frac{\partial \hat{u}}{\partial a_i} , F(\hat{u})\right \rangle_H =\left\langle \phi_i , F(\hat{u})\right \rangle_H. 
\end{equation}

Since the set of modes $\{\phi_i\}_{i=1}^N$ are linearly independent, the metric tensor $\vc M$ is symmetric positive-definite.  Moreover, if the modes are orthonormal then the metric tensor is the $N\times N$ identity matrix. In either case, the linear system of equations \eqref{eq:rons_ode} has a unique solution, $\dot{\vc a} = \vc M^{-1}\vc f(\vc a)$. 

\begin{figure}
\centering 
\includegraphics[width = .75\textwidth]{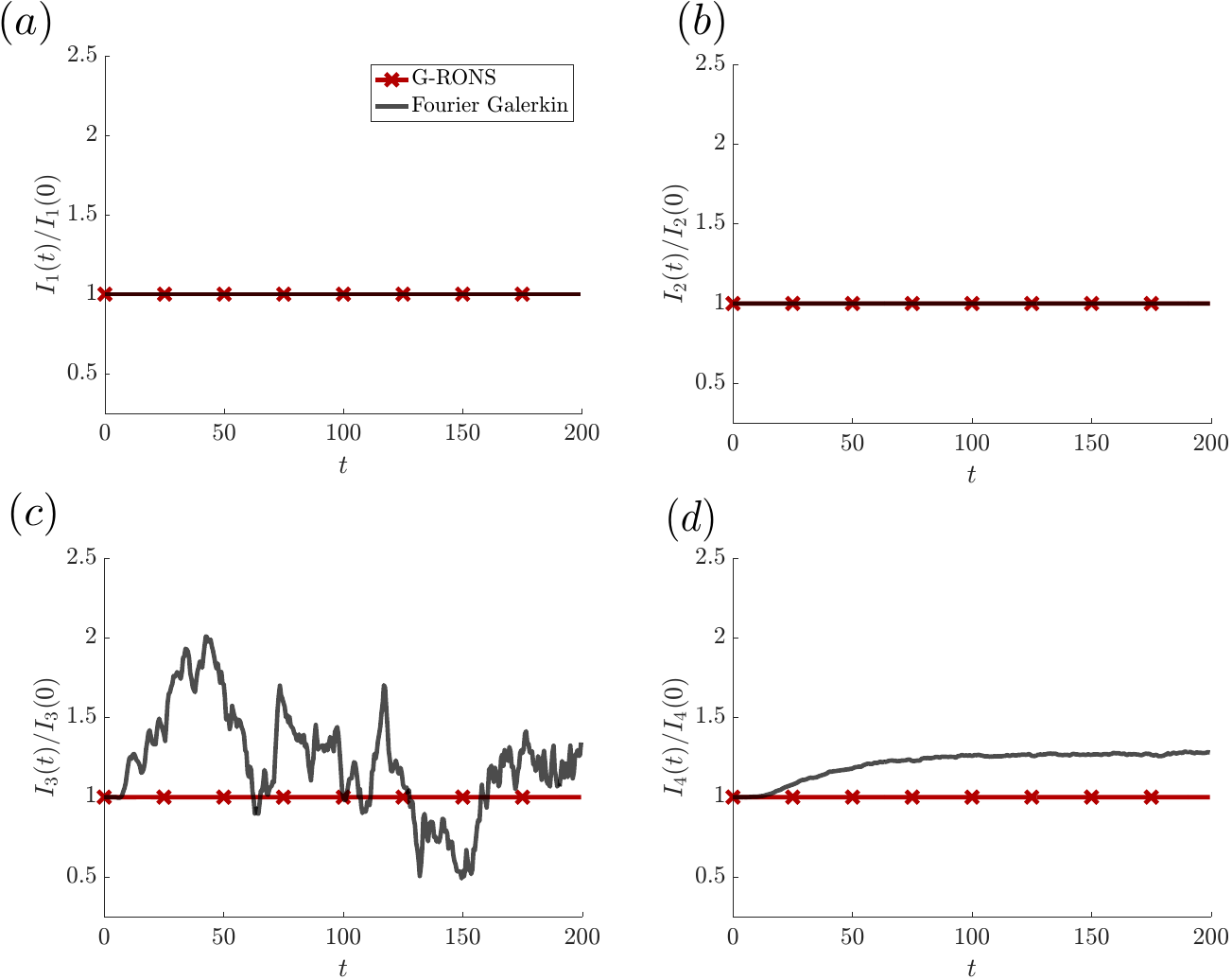}
\caption{\label{fig:euler}First four polynomial Casimir invariants for Euler's equation: (a) $I_1$, (b) $I_2$, (c) $I_3$, and (d) $I_4$. In each panel, the red line with the $\times$ markers corresponds to Galerkin RONS whereas the solid black curve denotes a standard Fourier psuedospectral DNS with $512 \times 512$ mdoes.}
\end{figure}
Once the ODE~\eqref{eq:rons_ode} is solved numerically, the approximate solution $\hat u(\vc x,\vc a(t))$ can be evaluated. However, there is no guarantee that this approximate solution
respects the conserved quantities~\eqref{eq:consQuantities}. More precisely, there is generally no guarantee that the quantities $I_k[\hat u(\cdot,\vc a(t))]$ are constant in time. Even in very high-order truncations, where $N$ is very large, the conserved quantities are not guaranteed to be respected. For instance, consider the two-dimensional Euler equation for ideal fluids with periodic boundary conditions, $\partial_t \omega+\vc u\cdot\nabla\omega=0$, where $\omega(\vc x,t)$ denotes the vorticity field. Euler's equation admits the Casimir invariants, $I_k = \int \omega^k(\vc x,t) \id\vc x$ for any integer $k\geq 1$. Let's consider a Fourier spectral approximation $\hat \omega$ of the solution where $\phi_i$ are the Fourier modes. In this approximation, only $I_1$ and $I_2$ remain conserved, however, all higher-order Casimir invariants $I_k$ with $k\geq 3$ are violated~\cite{Bowman2013}. This is shown in figure~\ref{fig:euler} for a Fourier pseudo-spectral solution of the Euler equation with $512\times 512$ modes over the time interval $t\in[0,200]$.

In the framework of RONS, it is straightforward to ensure that the conserved quantities are respected by the approximate solution. To this end, we add the conserved quantities as constraints to the optimization problem~\eqref{eq:rons_opt}. More specifically, we seek the solution to the constrained optimization problem,
\begin{subequations}\label{eq:rons_const_opt}
\begin{equation}
 \min_{\dot{\vc a}\in\R^N} \mathcal{J}(\vc a, \dot{\vc a}),
 \end{equation}
 \begin{equation}
\mbox{s.t.}\quad I_k(\vc a(t)) = const. \quad \text{for }  k = 1,2,\ldots , m.
  \end{equation}
\end{subequations}

The solution to this optimization problem can be written explicitly. Here, we recall the resulting equations and refer to~\cite{Anderson2022a} for a detailed derivation. 
The constraints in the optimization problem are enforced by introducing the Lagrange multipliers $\boldsymbol{\lambda}=(\lambda_1,\cdots,\lambda_m)^\top$ which solve the \textit{constraint equation},
\begin{equation} \label{eq:rons_lagrange}
\vc C(\vc a) \boldsymbol{\lambda} = \vc b(\vc a). 
\end{equation}
Here, the \textit{constraint matrix} $\vc C$ has entries 
\begin{equation}
C_{ij} =\left \langle \nabla_{\vc a} I_i,  \vc M^{-1} \nabla_{\vc a} I_j \right \rangle,
\end{equation}
and the components of  $\vc b$ are given by
\begin{equation}
b_i = \left \langle \nabla_{\vc a} I_i , \vc M^{-1} \vc f\right \rangle,
\end{equation}
where $\langle \cdot, \cdot \rangle$ denotes the usual Euclidean inner product on $\mathbb{R}^N$. 
\begin{rem} \label{rem:rons_grads}
The conserved quantities are typically defined by an integral over the spatial domain,
\begin{equation}
I(u) = \int_{\Omega}  g(u)\id \vc x,
\end{equation}
for some map $g$. After substituting the approximate solution $\hat u(\vc x,\vc a(t))$ into this equation and integrating over $\Omega$, the result is a function of $\vc a$.
Using a slight abuse of notation, we write $I(\vc a)$ for $I(\hat u(\cdot,\vc a))$. In the language of differential geometry, $I(\hat u(\cdot,\vc a))$ is the pullback of $I$ under the map $\hat u$.

%For Galerkin truncations, the derivatives needed for RONS can easily be found with a simple application of the chain rule: 
%\begin{equation}
%\frac{\partial }{\partial a_k} I(\vc a) = \int_\Omega g^\prime(\hat u(\vc x,\vc a)) \phi_k (\vc x)\id \vc x. 
%\end{equation}
\end{rem}
The constraint matrix $\vc C(\vc a)$ is symmetric positive-definite provided that the gradients $\left\{\nabla_{\vc a} I_k(\vc a)\right\}_{k=1}^m$ are linearly independent, so the solution to the constraint equation \eqref{eq:rons_lagrange} exists and is unique. With this set-up, the minimizer of \eqref{eq:rons_const_opt} must solve the ODEs,
\begin{equation}\label{eq:rons_eq}
\vc M \dot{\vc a} = \vc f (\vc a)- \sum_{k=1}^m \lambda_k \nabla_{\vc a} I_k(\vc a).
\end{equation}
We refer to Eq. \eqref{eq:rons_eq} as the Galerkin RONS equation, or G-RONS for short. The summation term involving the Lagrange multipliers $\lambda_k$ ensure that the corresponding approximate solution $\hat u(\vc x,\vc a(t))$
respects the conservation of the quantities $I_k$. We solve the RONS equation numerically using standard time discretization schemes such as Runge--Kutta methods. At each time step, the linear constraint equation~\eqref{eq:rons_lagrange} needs to be solved to obtain the Lagrange multipliers. Since the number of conserved quantities $m$ is often small, solving the constraint equation does not constitute a significant computational cost.

Comparing the usual Galerkin equation~\eqref{eq:rons_ode} and the Galerkin RONS equation~\eqref{eq:rons_eq}, we see that they only differ by the summation term involving the Lagrange multipliers. Therefore, it is quite straightforward to modify existing Galerkin code to enforce the conservation of the first integrals $I_k$.
For instance, figure~\ref{fig:euler} shows that Galerkin RONS conserves the higher-order Casimir invariants of the Euler equation.

\section{Extensions of RONS} \label{sec:extensions}
In this section, we develop two extensions of RONS. In section \ref{sec:RONS_systems}, we derive the Galerkin RONS equation for systems of PDEs. In section \ref{sec:fv_rons}, we apply RONS in conjunction with finite volume discretization to ensure that the discretized equations respect all conserved quantities of the PDE. 
\subsection{RONS for Systems of PDEs} \label{sec:RONS_systems}
The RONS theory can be generalized to systems of PDEs. Consider a system of $p$ PDEs in the form,
\begin{equation}\label{eq:vec_gen_pde}
\vc u_t  = \vc F (\vc u),
\end{equation} 
where $\vc u:\mathbb{R}^n\times\mathbb{R}^+\rightarrow \mathbb{R}^p$ has components $u_1,\ldots, u_p$, and $\vc F$ is a vector-valued differential operator with components $F_1, \ldots,F_p$. We assume that each component of the solution belongs to a Hilbert space $H_i$, i.e., $u_i\in H_i$.
As before, we assume that the PDE \eqref{eq:vec_gen_pde} has $m$ conserved quantities,
\begin{equation}
I_k(\vc u(\cdot, t)) = \mbox{const.}, \quad \ k=1,2,\ldots, m,
\end{equation}
for all times $t\geq 0$.

Consider an approximate solution $\hat{\vc u} = (\hat u_1,\cdots,\hat u_p)$ which depends on space $\vc x$ and time-dependent parameters $\vc a(t)$. 
More specifically, we write each component of the approximate solution as a linear combination of prescribed modes, 
\begin{equation}
	\hat{u}_k (\vc x,\vc a_k(t)) = \sum_{i=1}^{N_k}a_k^i(t)\phi_k^i(\vc x),
\end{equation} 
where $\{\phi^i_k(\vc x)\}_{i=1}^{N_k}$ are the modes expanding the component $\hat u_k$ of the approximate solution and $\{a^i_k(t)\}_{i=1}^{N_k}$ are the corresponding time-dependent parameters. We define the vector $\vc a_k = (a^1_k,\cdots,a^{N_k}_k)$ containing all parameters of the $k$-th component of $\hat{\vc u}$.
Note that for generality, we assumed that each component $\hat u_k$ is expanded by a different set of modes $\phi^i_k$. In the special case where all components are written in the same basis, the subscript $k$ can be omitted.

Concatenating all vectors $\vc a_k$ together, we define $\vc a = (\vc a_1,\cdots,\vc a_p)^\top\in \mathbb R^N$ where $N = N_1 + \ldots + N_p$.
Therefore, the approximate solution $\hat{\vc u}(\vc x,\vc a(t))$ depends on $N$ time-dependent parameters.
Similar to the scalar case, we define the cost functional,
\begin{equation}
\mathcal{J}(\vc a,\dot{\vc a}) = \frac{1}{2}\| \hat{\vc u}_t - \vc F(\hat{\vc u})\|^2_H,
\label{eq:systemPDE_J1}
\end{equation}
which measures the discrepancy between the rate of change of the approximate solution and the rate dictated by the right-hand side of the PDE.
Here the Hilbert space $H$ is defined by the direct sum $H =\bigoplus_{k=1}^{p} H_k$, where $u_k \in H_k$. The norm on $H$ is defined by $\| \vc u \|_H^2 := \sum_{k=1}^{p} \|u_k\|_{H_k}^2$. Therefore, we can expand the cost function and write
\begin{equation}
\mathcal{J}(\vc a,\dot{\vc a}) = \frac{1}{2}\sum_{k=1}^p \left\|\pard{\hat{u}_k}{t} - F_k(\hat{\vc u}) \right\|_{H_k}^2. 
\label{eq:systemPDE_J2}
\end{equation}
In order to express the cost function $\mathcal J$ more explicitly, we define the matrices $\vc M_k\in\mathbb R^{N_k\times N_k}$ and the vectors $\vc f_k(\vc a)\in\mathbb R^{N_k}$ for $k = 1,2,\ldots,p$. Denoting the components of $\vc M_k$ and $\vc f_k$ by $[\vc M_k]_{ij}$ and $[\vc f_k]_j$, respectively, we define,
\begin{equation} 
 [\vc M_k]_{ij} = \left\langle\phi^i_k,\phi^j_k \right\rangle_{H_k}, \quad [\vc f_k]_i =\left \langle \phi^i_k , F_k(\hat{u})\right \rangle_{H_k}, 
\end{equation}
where $ i,j\in\{1,2,\cdots,N_k\}$. Noting that,
\begin{equation}
\partial\hat u_k/\partial t= \sum_i \dot a_k^i(t)\phi_k^i(x), 
\end{equation}
it is straightforward to show that the cost functional~\eqref{eq:systemPDE_J2} can be written equivalently as
\begin{equation}
\begin{split}
\mathcal{J}(\vc a, \dot{\vc a}) =&  \sum_{k=1}^p \frac{1}{2}\left\langle\dot{\vc a}_k, \vc M_k \dot{\vc a}_k\right\rangle \\
& - \left\langle \dot{\vc a}_k ,\vc f_k \right\rangle + \frac{1}{2}\left\langle F_k(\hat{\vc u}), F_k(\hat{\vc u}) \right\rangle_{H_k}.
\end{split}
\label{eq:systemPDE_J3}
\end{equation}
Recall that $\langle\cdot,\cdot\rangle$ denotes the Euclidean inner product, whereas $\langle\cdot,\cdot\rangle_{H_k}$ denotes the functional inner product on the Hilbert space $H_k$.
Our final step is to write the cost functional~\eqref{eq:systemPDE_J3} more compactly by defining the block diagonal matrix $\vc M = \text{diag}\left(\vc M_1, \vc M_2, \ldots, \vc M_p \right)$, and the vector field $\vc f:\mathbb R^N\to\mathbb R^N$ such that
\begin{equation}
\vc M = \begin{bmatrix}
\vc M_1 & & &  \\
& \vc M_2 & & \\ 
& & \ddots \\
& & & \vc M_p 
\end{bmatrix}, \quad 
\vc f(\vc a) = \begin{bmatrix}
\vc f_1 (\vc a) \\
\vc f_2(\vc a)\\ 
\vdots \\
\vc f_p(\vc a) 
\end{bmatrix}.
\label{eq:systemPDE_Mf}
\end{equation}
With these expressions, the cost function takes its final form,
\begin{equation}
\begin{split}
\mathcal{J}(\vc a, \dot{\vc a}) & = \frac{1}{2}\left\langle \dot{\vc a}, \vc M\dot{\vc a}\right\rangle - \left\langle \dot{\vc a}, \vc f \right\rangle \\ 
 &+ \frac{1}{2}\sum_{k=1}^p \left\langle F_k(\hat{\vc u}),F_k(\hat{\vc u}) \right\rangle_{H_k},
\end{split}
\end{equation}
which is similar in form to the scalar-valued case reviewed in section~\ref{sec:math_prelim}.

We minimize this cost function instantaneously. To enforce the conserved quantities of the PDE, we follow the same procedure as in the scalar-valued case by adding them as constraints to the instantaneous optimization problem. More specifically, we seek the solution of the constrained optimization problem,
\begin{align}
	\begin{split}
		\min_{\dot{\vc a}\in\R^N} \  & \mathcal{J}(\vc a, \dot{\vc a}), \\
		\text{s.t.} \ & I_k(\vc a(t)) = \mbox{const.} \quad \text{for }  k = 1,2,\ldots , m.
	\end{split}
\end{align}
The solution to this constrained optimization problem satisfies the ODEs,
\begin{equation}\label{eq:rons_systemPDE}
	\vc M \dot{\vc a} = \vc f(\vc a)- \sum_{k=1}^m \lambda _ k  \nabla_{\vc a} I_k (\vc a),
\end{equation}
with $\vc M$ and $\vc f$ defined in~\eqref{eq:systemPDE_Mf}. The proof is similar to the scalar-valued case and therefore is omitted here for brevity.

Note that we massaged the cost function $\mathcal J$ for the system of PDEs so that it has a similar form to the scalar case. Therefore, it is not surprising that the resulting equations are identical in form to equation~\eqref{eq:rons_eq}. However, there are notable differences. In particular, here $\vc M$ is the block-diagonal matrix defined in~\eqref{eq:systemPDE_Mf} and similarly $\vc f$ is the vector field obtained by concatenating the vectors $\vc f_k$.
Since the modes $\{\phi_k^i\}_{i=1}^{N_k}$ are linearly independent, each block matrix  $\vc M_k$ is symmetric positive definite. Hence, the block-diagonal matrix $\vc M$ is also symmetric positive definite and therefore invertible.

The Lagrange multiplier $\bs \lambda=(\lambda_1,\cdots,\lambda_m)^\top$ solves the linear constraint equation $\vc C(\vc a) \bs \lambda  = \vc b(\vc a)$, as in section \ref{sec:math_prelim}. The components of the constraint matrix $\vc C$ are given by $C_{ij} = \langle\nabla_{\vc a} I_i, \vc M^{-1} \nabla_{\vc a} I_j \rangle$, where the gradients of the conserved quantities $I_k$ are now defined by
\begin{equation}
\nabla_{\vc a} I _ k(\vc a)= \begin{bmatrix}
\nabla_{\vc a_1 } I_k (\vc a)\\ 
\nabla_{\vc a_2 } I_k (\vc a)\\ 
\vdots \\
\nabla_{\vc a_p } I_k (\vc a)\\ 
\end{bmatrix}.
\end{equation}
The vector $\vc b$ is defined in the same way as in section \ref{sec:math_prelim}, i.e., $b_i = \langle \nabla_{\vc a} I_i, \vc M^{-1}\vc f \rangle$. Assuming that the set of all gradients $\{\nabla_{\vc a} I_k\}_{k=1}^m$ is linearly independent, then the constraint matrix $\vc C$ is symmetric positive-definite. Thus the solution to the constraint equation exists and is unique.

We conclude this section by noting that the above theory can be easily extended to the case where the components of the approximate solution $\hat{u}_k(\vc x, \vc a_k(t))$ depend nonlinearly on their corresponding parameters $\vc a_k$. In this case, one needs to define $\vc f_k$ and $\vc M_k$ differently as outlined in \cite{Anderson2022a}. But the remainder of the derivation remains unchanged if we allow the nonlinear dependence of $\hat{u}_k$ on the parameters $\vc a_k$.

\subsection{Finite Volume RONS}\label{sec:fv_rons}
In this section, we consider finite volume methods for direct numerical simulation of PDEs modeling conservation laws. 
Consider a general conservation law of the form,
\begin{equation}\label{eq:cons_law}
	u_t + \nabla \cdot \vc G(u) = 0,
\end{equation}
where $u(\vc x,t)\in\mathbb R$ is the solution to the PDE and $\vc G$ is a potentially nonlinear differential operator. Note that, for simplicity, we again assume that the solution is scalar-valued. Nonetheless, our results can be extended to systems of conservation laws using the method introduced in section~\ref{sec:RONS_systems}.

Under certain conditions, equation~\eqref{eq:cons_law} conserves the spatial integral of the state variable $u$. For example, if the boundary conditions are periodic, or if $\vc G(u)$ is orthogonal to the unit normal vector $\vc n$ of the boundary, then
\begin{equation} 
\begin{split}
\pard{}{t}\int_{\Omega} u(\vc x,t)\id \vc x &  = - \int_{\Omega} \nabla \cdot \vc G(u) \id\vc x \\
& = -\int_{\partial \Omega} \vc G(u)\cdot \vc n \id S = 0,
\end{split}
\end{equation}
which implies the conservation of the state variable $\int_{\Omega} u(\vc x, t) \id\vc x$.
Certain conservation law PDEs admit additional conserved quantities which are not so apparent. For instance, the shallow water equation, discussed in section~\ref{sec:swe},
conserves the total energy. This additional conserved quantity can be inferred from the Hamiltonian structure of shallow water equation~\cite{Camassa1994}, but it is not apparent from its conservative formulation~\eqref{eq:cons_law}.

Finite volume schemes are designed to conserve the state variable \cite{Chertock2015,Carlberg2018,Kurganov2018}. However, there is no guarantee that the numerical solution respects the conservation of additional invariants, e.g., the total energy for the shallow water equation. As reviewed in section~\ref{sec:rel_work}, there are structure-preserving schemes which enforce the conservation of additional conserved quantities~\cite{Fjordholm2011,Tadmor2008,Mishra2011}. However, these schemes are typically derived for a particular equation and are not necessarily applicable to different systems. Here, our goal is to apply the RONS framework to finite volume schemes in order to ensure that the discretized equations respect the conservation of any additional invariants. 
As we show below, our method is broadly applicable to any finite volume scheme and is not limited to a particular PDE.

%Many problems in fluid mechanics \cite{Mishra2011} and ocean and atmospheric sciences \cite{Chertock2015,Fjordholm2011}, among others, can be solved numerically using finite volume.  

First, we review the general formulation of the finite volume method.
As opposed to working with the strong formulation of the conservation law \eqref{eq:cons_law}, finite volume methods discretize its integral formulation. More specifically, let $\Omega\subseteq \R^n$ be the spatial domain of the PDE, and consider the set of $N$ time-independent and disjoint control volumes $\Omega_i$, such that $\overline{\bigcup_{i=1}^N \Omega_i} = \Omega$. Consider the integral of \eqref{eq:cons_law} over the $i$-th control volume, 
\begin{equation}\label{eq:fv_1}
\int_{\Omega_i} u_t \id \vc x + \int_{\Omega_i} \nabla \cdot \left(\vc G(u)\right) \id \vc x = 0.
\end{equation}
The cell average, over this control volume, is defined by
\begin{equation}
U_i(t) = \frac{1}{|\Omega_i|}\int_{\Omega_i} u(\vc x,t) \id\vc x.
\label{eq:cellAvg}
\end{equation}
where $|\Omega_i|$ denotes the volume of the $i$-th cell. Since the cells are independent of time, the integral of the conservation law on $\Omega_i$ can  be expressed in terms of the cell averages,
\begin{equation}
\begin{split}
\frac{\id}{\id t}U_i &= -\frac{1}{|\Omega_i|} \int_{\Omega_i}\nabla \cdot \vc G(u) \id \vc x \\ 
& =- \frac{1}{|\Omega_i|}\int_{\partial\Omega_i}  \vc G(u)  \cdot \vc n \id S,
\end{split}
\end{equation}
where we applied the divergence theorem with $\vc n$ being the outward-facing unit normal vector of $\Omega_i$. Therefore, the time derivative of the cell average $U_i$ can be written as the flux of $\vc G(u)$ through the boundary of $\Omega_i$. Summing over all cells, the flux through each internal boundary $\partial\Omega_i$ appears twice with opposite signs, since the unit normals point in opposite directions. This fact ensures the conservation of the state variable $\int_\Omega u\id\vc x = \sum_i |\Omega_i|U_i$.

\begin{comment}
Next consider the sum of the time derivatives of the cell averages $\sum_{i=1}^NU_i$; two neighboring cells $U_i$ and $U_j$ will have unit normals $\vc n$ in opposing directions. Thus the sum of these two fluxes must be equal to zero. For a more concrete example, consider a one dimensional system in the form \eqref{eq:cons_law}. In this case, each  control volume is a sub-interval of $\Omega$. Therefore the flux across the boundary of $\Omega_i$ is the flux across the left and right endpoints of $\Omega_i$ or $\mathcal{G}_i(\vc U)$ and $\mathcal{G}_{i+1}(\vc U)$ respectively. Thus the time derivative of the $i$th cell average is 
\begin{equation}
\frac{\id}{\id t} U_i = -\mathcal{G}_{i+1}(\vc U)+ \mathcal{G}_{i}(\vc U),
\end{equation}
where $\vc U$ is the vector whose $i$th component is $U_i$. If we consider the sum of the time derivatives of two adjacent cells, we find that 
\begin{equation}
\frac{\id}{\id t}U_i +\frac{\id}{\id t}U_{i+1}  = -\mathcal{G}_{i+1}(\vc U)+ \mathcal{G}_{i}(\vc U) -\mathcal{G}_{i+2}(\vc U)+ \mathcal{G}_{i+1}(\vc U) = \mathcal{G}_{i}(\vc U) -\mathcal{G}_{i+2}(\vc U).
\end{equation}
Taking the sum of $\frac{\id}{\id t}U_i$ over all $N$ subdomains, the only terms that survive are the values of the fluxes across the boundary of the entire domain $\Omega$. Under  certain conditions such as periodic boundary conditions, the difference between the fluxes on the boundary of $\Omega$ is zero. In this case, the state variable $u$ is conserved on the discrete level. 
\end{comment}

In the finite volume method, the flux through the boundary $\partial \Omega_i$ is approximated using the values of the $i$-th cell average and those of the adjacent cells. 
Different finite volume schemes correspond to different methods for estimating these fluxes. In general, the approximation of the average flux through the $i$-th cell can be written as a function $\mathcal F_i(\vc U)$ where $\vc U = (U_1,\cdots, U_N)^\top$ denotes the vector containing all cell averages. The corresponding finite volume scheme is then given by
\begin{equation}
\frac{\id }{\id t}U_i = \mathcal{F}_{i}(\vc U),\quad i=1,2,\cdots,N.
\label{eq:discFV}
\end{equation}

Now we turn our attention to finite volume RONS. 
%Our starting point is the finite volume discretization~\eqref{eq:discFV}.
For the original PDE~\eqref{eq:cons_law}, we consider an approximate solution,
\begin{equation}\label{eq:fvRONS_uhat}
	\hat{u}(\vc x,t) = \sum_{i=1}^N U_i(t) \phi_i(\vc x),
\end{equation} 
where $U_i$ are the cell averages~\eqref{eq:cellAvg} and $\phi_i$ is the indicator function on $\Omega_i$, i.e.,
\begin{equation}
\phi_i(\vc x) = \begin{cases}
1 & \vc x \in \Omega_i, \\
0 & \vc x \notin \Omega_i.
\end{cases}
\end{equation}
Being an approximate solution, $\hat u$ does not exactly solve the PDE~\eqref{eq:cons_law}.
Next we define $F(\hat u) = \sum_{i=1}^N \mathcal{F}_i(\vc U)\phi_i(\vc x)$. Using the finite volume discretization~\eqref{eq:discFV}, the approximate solution $\hat u$ solves the approximating PDE,
\begin{equation}\label{eq:approx_PDE}
	\partial_t \hat u = F(\hat u).
\end{equation}
Although this PDE differs from the original PDE~\eqref{eq:cons_law}, it admits the same finite volume discretization. In fact, taking the spatial integral of~\eqref{eq:approx_PDE} over the control volume $\Omega_i$, we recover equation~\eqref{eq:discFV} exactly. 

We apply RONS to the approximating PDE~\eqref{eq:approx_PDE}. More precisely, as in section~\ref{sec:math_prelim}, we solve the constrained optimization problem,
\begin{subequations}\label{eq:constOpt_FV}
\begin{equation}
	\min_{\dot{\vc U}\in\R^N} \mathcal J(\vc U,\dot{\vc U}) : = \frac{1}{2}\|\hat{u}_t - F(\hat u) \|_{L^2}^2,
\end{equation}
\begin{equation}
		\mbox{s.t.}\quad I_k(\vc U(t)) = \mbox{const.} \quad \text{for }  k = 1,\ldots , m,
\end{equation}
\end{subequations}
where $I_k$ denote all conserved quantities of the original PDE~\eqref{eq:cons_law}. These include the state variable $\int_\Omega \hat u(\vc x,t)\id \vc x = \sum_i |\Omega_i|U_i(t)$ 
and any additional conserved quantities. Using the explicit form of the cost function, 
\begin{equation}
\begin{split}
	\mathcal{J}(\vc U, \dot{\vc U}) = \frac{1}{2}\sum_{i=1}^N  |\Omega_i|  \biggl \{  & \dot{U}_i^2   - 2\dot{U}_i\mathcal{F}_i(\vc U)  \\
	& + [\mathcal{F}_i(\vc U)]^2\biggr\}, 
\end{split}
\end{equation}
and applying RONS, the solution to the constrained optimization problem~\eqref{eq:constOpt_FV} is given by
\begin{equation}\label{eq:rons_eq_fv}
	\dot{\vc U} = \pmb{\mathcal{F}}(\vc U) - \sum_{k=1}^m\lambda_k \vc M^{-1}\nabla_{\vc U} I_k(\vc U),
\end{equation}
where the vector field $\pmb{\mathcal F} = (\mathcal F_1,\cdots,\mathcal F_N)^\top$ contains the fluxes $\mathcal F_i$ through the 
control volumes $\Omega_i$.
Here, the metric tensor $\vc M$ is diagonal and positive definite, $M_{ij} = |\Omega_i|\delta_{ij}$, and therefore its inverse is computable at insignificant additional computational cost.

As in section~\ref{sec:math_prelim}, the Lagrange multipliers $\pmb\lambda = (\lambda_1,\allowbreak,\cdots,\lambda_m)^\top$ satisfy the linear system $\vc C \bs \lambda = \vc b$, where
the constraint matrix $\vc C$ and the vector $\vc b$ are given by  $C_{ij} = \langle\nabla_{\vc U} I_i,\vc M^{-1}\nabla_{\vc U}I_j \rangle$ and $b_i = \langle\nabla_{\vc U} I_i,\pmb{\mathcal{F}} \rangle$, respectively. If no conserved quantities are enforced, i.e., $\pmb\lambda=\vc 0$, then finite volume RONS coincides with the classical finite volume discretization~\eqref{eq:discFV}.
We refer to~\eqref{eq:rons_eq_fv} as the finite volume RONS equation, or FV-RONS for short. Note that the form of this equation is slightly different from the Galerkin RONS equation~\eqref{eq:rons_eq} since, for FV-RONS, the equations simplify due to the diagonal metric tensor $\vc M$.

Let $I_1$ denote the integral of the state variable, i.e., 
\begin{equation}
	I_1(\vc U)  = \int_\Omega \hat u\id \vc x = \sum_i | \Omega_i| U_i.
\end{equation}
The gradient of this conserved quantity is related to the cell volumes, 
\begin{equation}
	\pard{I_1}{U_i} = |\Omega_i|.
\end{equation}
The gradient of any additional conserved quantities depends on the form of that conserved quantity.

In section~\ref{sec:swe}, we present numerical results showing that finite volume RONS conserves both the state variable and the Hamiltonian (i.e., total energy) of the system, whereas the classical finite volume methods only conserve the state variable.

\section{Numerical Results}\label{sec:num_results}
Here we present two numerical examples: the shallow water equation (SWE) and the nonlinear Schr\"odinger (NLS) equation. For the shallow water equation, we use finite volume RONS to enforce the conservation of total energy in a central-upwind finite volume scheme. As such, our goal in section~\ref{sec:swe} is not reduced-order modeling; rather the solutions are full-order and we compare finite volume RONS against a classical finite volume method. In section~\ref{sec:nls}, we study reduced-order models for the nonlinear Schr\"{o}dinger equation where we demonstrate the effects of enforcing conserved quantities in the reduced model.

\subsection{Shallow Water Equation}\label{sec:swe}
The shallow water equation, also referred to as the Saint-Venant system, is a system of PDEs that describes the evolution of surface waves when the depth of the fluid is much smaller than the dominant wavelength of the surface waves. SWE is used in modeling the behavior of tsunami's \cite{Kevlahan2019,Khan2022}, as well as other phenomenon such as the dam break problem in engineering \cite{Wang2020}. In one spatial dimension and in dimensionless variables, the conservative form of the shallow water equation is given by
\begin{align}\label{eq:swe}
\begin{split}
\eta_t + \left(\left(\eta+H(x)\right)v\right)_x &= 0, \\
v_t + \left(\frac{1}{2}v^2 + \gr\eta \right)_x  & = 0,
\end{split}
\end{align}
with the corresponding initial data $\eta(x,0)=\eta_0(x)$ and $v(x,0)=v_0(x)$. 
\begin{figure}
\centering
\includegraphics[width=.65\textwidth]{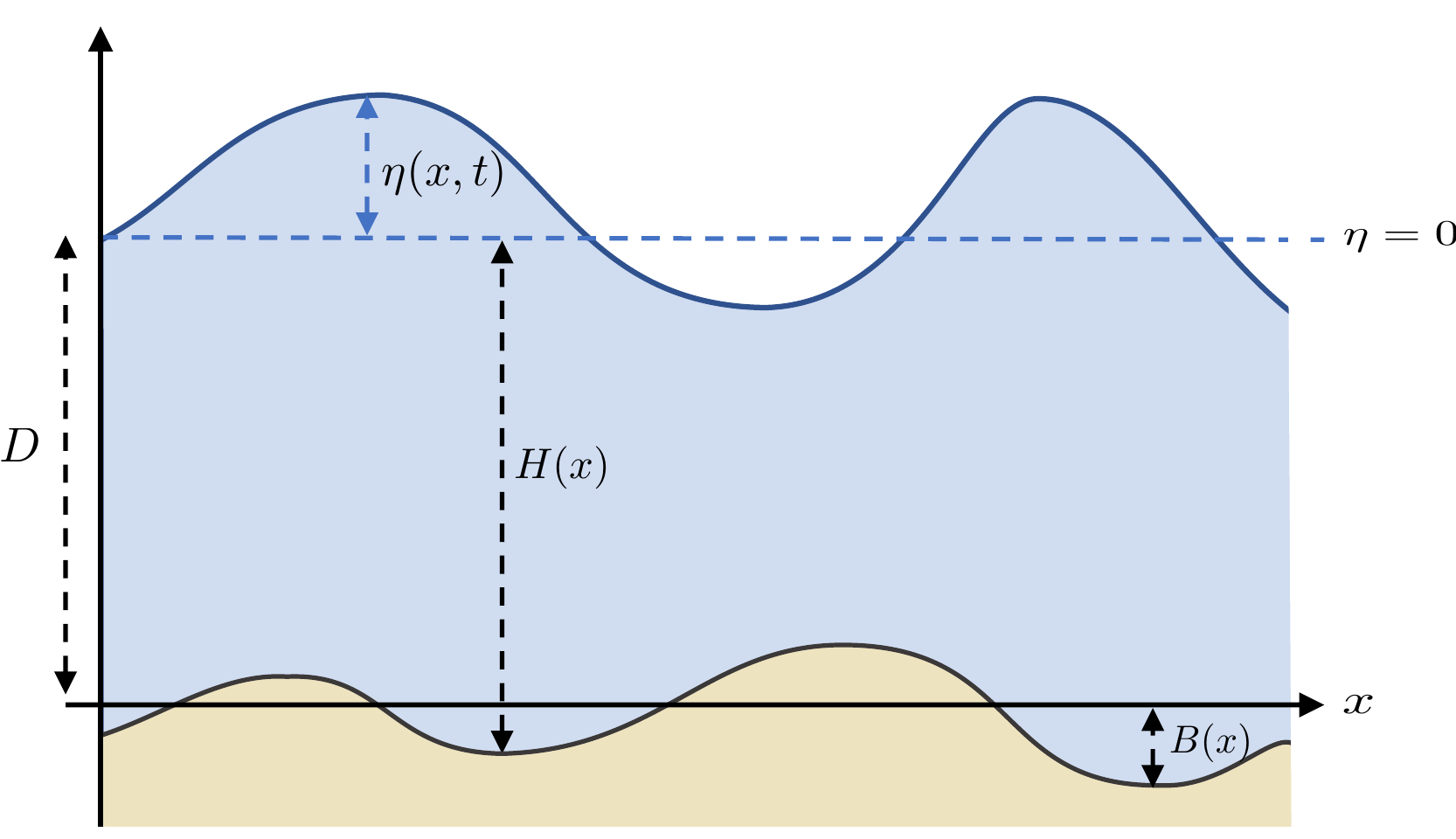}
\caption{\label{fig:swe_diagram} Schematic illustration of the set-up for the shallow water equation.}
\end{figure}
Here $v(x,t)$ is the depth-averaged horizontal velocity, $\gr$ is the acceleration due to gravity, and $\eta(x,t)$ is the surface elevation away from the surface at rest $\eta = 0$. As depicted in Figure \ref{fig:swe_diagram}, $H(x)$ is the distance from the fluid surface at rest to the bottom, and the function $B(x)$ defines the bottom topography away from the mean depth $D$ of the fluid. The depth $H(x)$, the bottom topography $B(x)$, and the average depth $D$ are related by $H(x) = D- B(x)$. Here, we assume periodic boundary conditions over a domain $x\in\Omega = [0,L]$.

Equation~\eqref{eq:swe} is written in dimensionless variables. If we denote the dimensional variables with a tilde, the dimensionless variables are defined by 
$x = \tilde{x}/\ell$, $t = \bar{v}\,\tilde{t}/\ell$,  $\eta = \tilde{\eta}/\ell$, $v =\tilde{v}/ \bar{v}$, 
where $\ell$ is the characteristic wavelength and $\bar{v}$ is the characteristic wave speed. As determined by the Intergovernmental Oceanographic Commission~\cite{tsunami_charvars}, we choose the characteristic wavelength $\ell = 2.13\times10^6$~m and speed $\bar{v} = 198$~m/s, which correspond to the average wavelength and speed of a tsunami wave in the ocean at depth $D = 4000$ m.  Lastly, the dimensionless gravitation acceleration is taken to be $\gr = 9.8\, \ell/\bar{v}^2$. 
 
The shallow water equation is a hyperbolic set of PDEs which is well-known to develop shocks and other discontinuities. However, under certain conditions, global strong solutions to SWE exist and are unique~\cite{Pelinovsky2017,Bai2022}. Here, we only consider such strong solutions which are known to conserve the following quantities,
 \begin{subequations}\label{eq:cq_swe}
 	\begin{equation} 
 		I_1(t) = \int_{\Omega}\eta(x,t) \id x, 
 		\label{eq:I_eta}
 	\end{equation}
\begin{equation} 
	I_2(t)= \int_{\Omega}v(x,t) \id x, 
	\label{eq:I_u}
\end{equation}
\begin{equation}
	I_3(t) = \frac{1}{2}\int_\Omega \left[(\eta + H)v^2 + \gr\eta^2\right] \id x.
	\label{eq:swe_energy}
\end{equation}
\end{subequations}
The first two quantities, $I_1$ and $I_2$, corresponding to the state variables, are trivially conserved. The third quantity $I_3$ is the total energy, or the Hamiltonian, of the system. A straightforward calculation shows that the total energy is also conserved along strong solutions of the shallow water equation~\cite{Camassa1994}.

In this section, our goal is to demonstrate that traditional finite volume schemes, which are designed to conserve the state variables, fail to conserve the total energy $I_3$. 
On the other hand, finite volume RONS, introduced in section~\ref{sec:fv_rons}, respects all these conserved quantities. Furthermore, we discuss the implications of this additional conserved quantity.

To this end, we first describe the traditional finite volume scheme used here. We use a second-order central-upwind finite volume scheme as developed in~\cite{Kurganov2007}. The spatial domain is $\Omega = [0,10]$ which we divide into $2^{10}$ cells of equal length. At each time step, we use a linear polynomial reconstruction to compute the fluxes with the nonlinear minmod flux limiter,
\begin{equation}
\begin{split}
\tilde{U_i}(x) &= U_i +  \frac{\Delta x}{2}\text{minmod}(\theta \frac{U_i-U_{i-1}}{\Delta x},\\
 &  \theta \frac{U_{i+1}-U_{i-1}}{2\Delta x}, \theta \frac{U_{i+1}-U_{i}}{\Delta x})(x-x_i),
\end{split}
\end{equation}
where $x_i$ is the center of the $i$-th control volume and $\Delta x$ is the size of each cell. We take the parameter $\theta = 1.2$ to ensure that the reconstructions are minimally oscillatory~\cite{Kurganov2007}. For the time integration, we use a third order strong-stability-preserving Runge-Kutta (SSP-RK3) scheme in accordance with Refs.~\cite{Chertock2015,Kurganov2018,Kurganov2007}. This scheme is described in detail in section 2.4.2 of Ref.~\cite{Gottlieb2011}. At each iteration of SSP-RK3, we choose the time step to satisfy the appropriate Courant--Friedrichs-Lewy (CFL) condition. More specifically, we choose the time step based on the maximum of the eigenvalues of the Jacobian of the map,
$$
\vc G(\eta,v) = \begin{bmatrix}
(\eta + H)v \\
\frac{1}{2}(v)^2 + g\eta
\end{bmatrix}.
$$
The eigenvalues are given by
\begin{equation}
\begin{split}
\lambda_1(x,t)& =  v(x,t) + \sqrt{g\left(\eta(x,t) + H(x)\right)},\\ 
 \lambda_2(x,t) &=  v(x,t) - \sqrt{g\left(\eta(x,t) + H(x)\right)},
\end{split}
\end{equation}
which should not be confused with the Lagrange multipliers in the RONS equation~\eqref{eq:rons_eq_fv}.
We set the time step $\Delta t_i$ at the $i$-th iteration of SSP-RK3 as 
\begin{equation}
\Delta t_i =\frac{\Delta x}{2\max\left\{\lambda_1^{\max}(t_{i-1}),\lambda_2^{\max}(t_{i-1})\right\}},
\end{equation}
where $\lambda_1^{max}(t_{i-1})= \max_x\{\lambda_1(x,t_{i-1})\}$ and $\lambda_2^{max}(t_{i-1}) = \max_x\{-\lambda_2(x,t_{i-1})\}$.

Finite volume RONS uses the same discretization with the important difference that the summation term in equation~\eqref{eq:rons_eq_fv} is added to the right-hand side of the spatial discretization, enforcing the conservation of the first integrals. Here, we have $m=3$ corresponding to three conserved quantities in~\eqref{eq:cq_swe}.

First, we test the finite volume (FV) method and FV-RONS on two benchmark examples: the lake at rest~\cite{Chertock2015} and the Gaussian pulse~\cite{Kevlahan2019,Khan2022}. The lake at rest is a steady state solution where both $\eta(x,t)$ and $v(x,t)$ are zero. Both finite volume and FV-RONS were able to capture this steady state with flat and non-flat bottom topographies (not shown here). For the Gaussian pulse, the initial conditions are given by
\begin{equation}\label{eq:gauss_pulse}
v_0(x) = 0,  \qquad \eta_0(x) = \frac{0.1}{\ell} e^{-\left(5(x-5)\right)^2}. 
\end{equation}
We integrate these initial conditions for ten time units. The initial amplitude of the wave is small enough that no shock waves form over this time interval, i.e., strong solutions exist.

\begin{figure*}
\centering
\includegraphics[width = \textwidth]{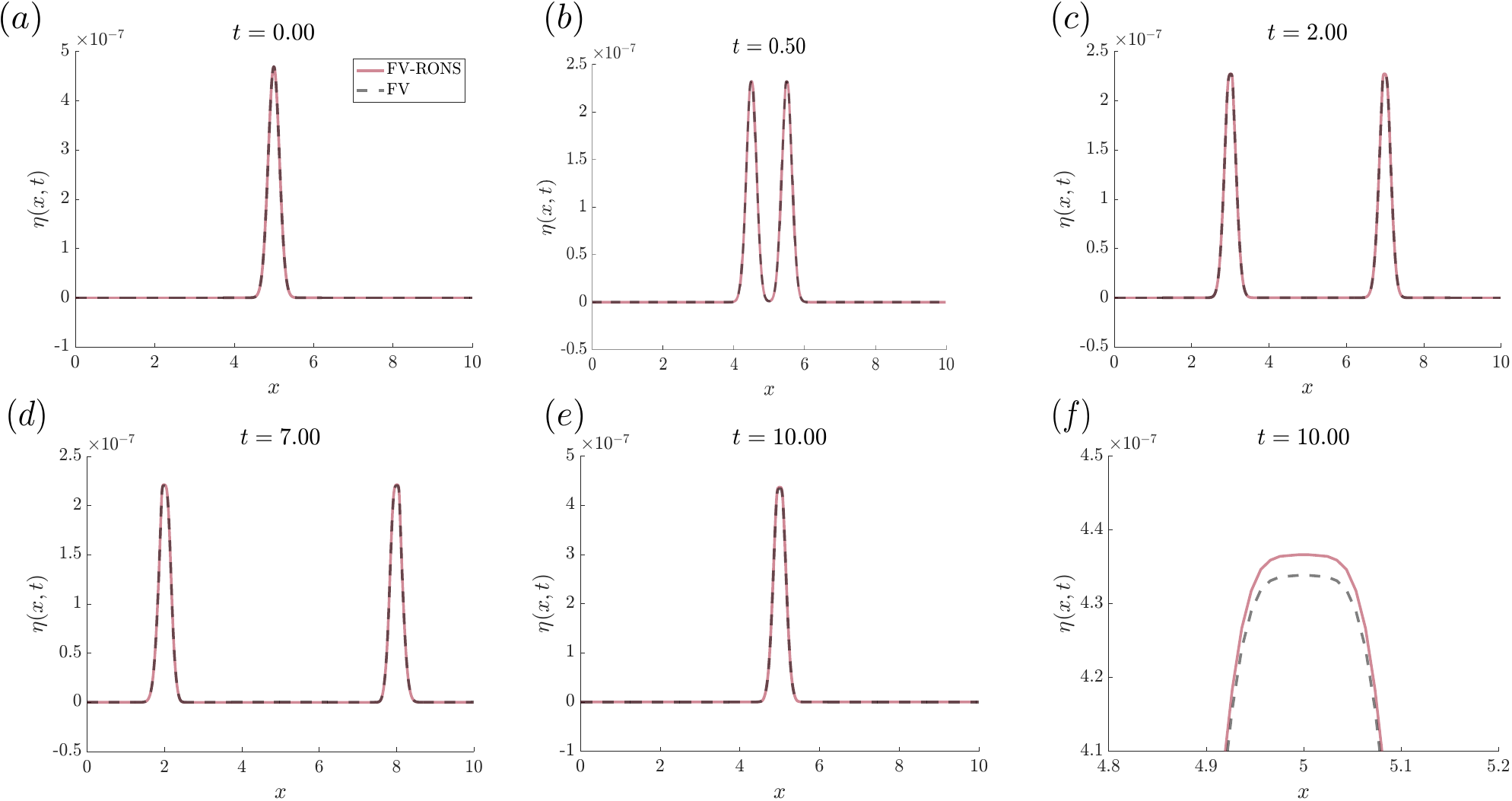}
\caption{Evolution of the Gaussian pulse initial condition, defined in Eq.~\eqref{eq:gauss_pulse} with the central-upwind finite volume (FV) scheme (dashed black line) and FV-RONS (solid red line) at times (a) $t = 0$, (b) $t = 0.5$, (c) $t = 2 $, (d) $t = 7$, (e) $t = 10$, (f) close-up of the pulses as they pass through one another at $t=10$.}
\label{fig:swe_gauss}
\end{figure*}
As shown in figure~\ref{fig:swe_gauss}, as the pulse evolves, we see similar qualitative behaviors for both FV and FV-RONS. Both solutions behave as expected for this initial condition: initially, the Gaussian splits into two smaller pulses that translate in opposite directions (fig.~\ref{fig:swe_gauss}(a)-(c)). They pass through the periodic boundary and then travel towards the center where they eventually merge with one another (fig.~\ref{fig:swe_gauss}(d) -(e)). An important difference between the two solutions is visible in Figure \ref{fig:swe_gauss}(f) which shows a close-up view of the solution at time $t=10$. The amplitude of the FV solution has decayed slightly more than the FV-RONS solution since FV does not conserve the total energy while FV-RONS does.

To further illustrate the differences between FV and FV-RONS, we turn our attention to a more oscillatory initial condition. Specifically, we define
\begin{equation}\label{eq:ic_eta}
\tilde{\eta}_0(x) = \cos(2\pi x)\left(\sum_{j=2}^5 \alpha_j\cos\left(2\pi j x + \varphi_j\right) \right),
\end{equation}
where the amplitudes $\alpha_j$ are sampled from a standard normal distribution, and the phases $\varphi_j$ are sampled from a uniform distribution on $[0,2\pi]$. To ensure that the amplitude of the initial surface elevation is small enough, allowing for strong solution to exist for a long time, we rescale $\tilde \eta_0$ and define the initial surface elevation,
\begin{equation}
%\eta_0(x) = \frac{A\tilde{\eta}_0(x)}{\ell \max_x\{\tilde{\eta_0}(x)\}}. 
\eta_0(x) = \frac{\tilde{\eta}_0(x)}{2 \ell \max_x\{\tilde{\eta_0}(x)\}}. 
\end{equation}
%Additionally, since tsunami amplitudes in deep water are roughly between the order of centimeters to meters \cite{tsunami_charvars}, we take the value of $A = 0.5$~m. 
For the initial velocity we set $v_0(x)= 0$. We evolve these initial conditions over the time interval $t\in[0,75]$.
\begin{figure*}
\centering
\includegraphics[width = \textwidth]{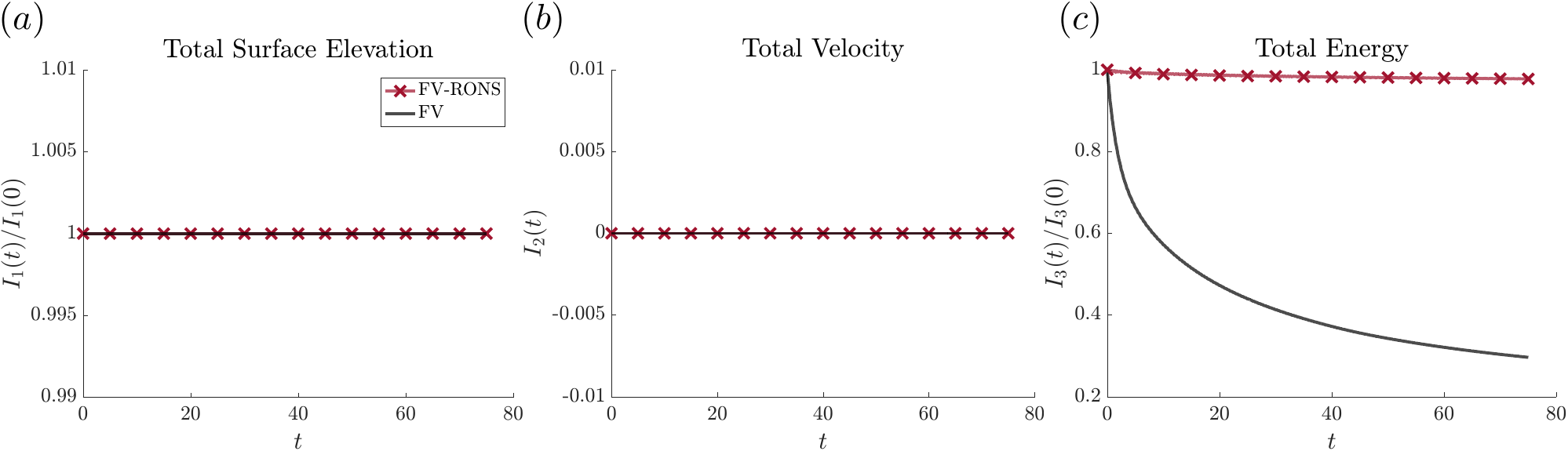}
\caption{Three conserved quantities defined in \eqref{eq:I_eta} - \eqref{eq:swe_energy} for finite volume (FV) and finite volume  
	RONS (FV-RONS). (a) The normalized total surface elevation $I_1$, (b) the total velocity $I_2$, and (c) the normalized total energy $I_3$. In each panel, the solid black line denotes FV while FV-RONS is the solid red line.}
\label{fig:fvrons_cons_q}
\end{figure*}

As shown in figure~\ref{fig:fvrons_cons_q}, both FV method and FV-RONS conserve the state invariants $I_1$ and $I_2$. However, the total energy $I_3$ decays significantly when using the FV method. In contrast, FV-RONS respects the conservation of the total energy. We emphasize that, as reviewed in section~\ref{sec:rel_work}, there exist energy-preserving finite volume methods. Our purpose here is to show that FV-RONS can preserve energy by minimally modifying existing FV methods as opposed to developing and implementing specialized FV schemes.
\begin{figure*}[htb]
\centering
\includegraphics[width = \textwidth]{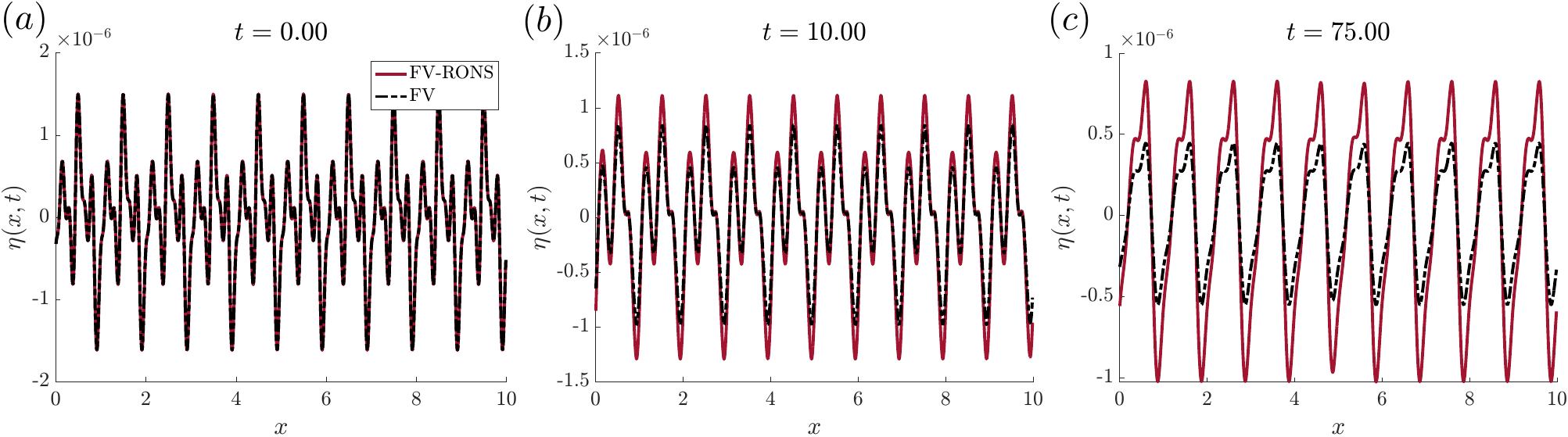}
\caption{Comparisons of the surface elevations $\eta(x,t)$ for FV (dashed black) and RONS (solid red) at different snapshots in time, (a) $t=0$, (b) $t=10$, and (c) $t=75$.}
\label{fig:fvrons_etas}
\end{figure*}

The energy decay of the FV method has a significant impact on the amplitude of the solutions. 
Figure~\ref{fig:fvrons_etas} compares the evolution of the numerical solutions at three time instances. Although the initial conditions are identical, significant differences begin to emerge quickly. 
At time $t=10$, although both methods predict the same wavelength, the amplitude of the FV solution has decayed noticeably more than FV-RONS. This behavior persists for later times; namely by time $t=75$, the amplitude of the FV solution is almost half of the amplitude of the FV-RONS solution.
This amplitude decay is reminiscent of the results shown in figure~\ref{fig:swe_gauss}(f); the decay is more significant when the waves are more oscillatory. 

Next, we compared the two numerical methods over a varied range of initial conditions. We sample $10^4$ different initial conditions from~\eqref{eq:ic_eta}, with randomly drawn amplitudes $\alpha_j$ and phases $\varphi_j$, and evolve them over the time interval $[0,75]$.  We compute the probability density function (PDF) of the maximum surface elevation $\max_{x\in \Omega} |\eta(x,t)|$. To this end, we discard the transients and only record the surface elevation for $t\in[25,75]$ at every $\Delta t = 0.1$ time units. 
\begin{figure}
\centering
\includegraphics[width=0.45\textwidth]{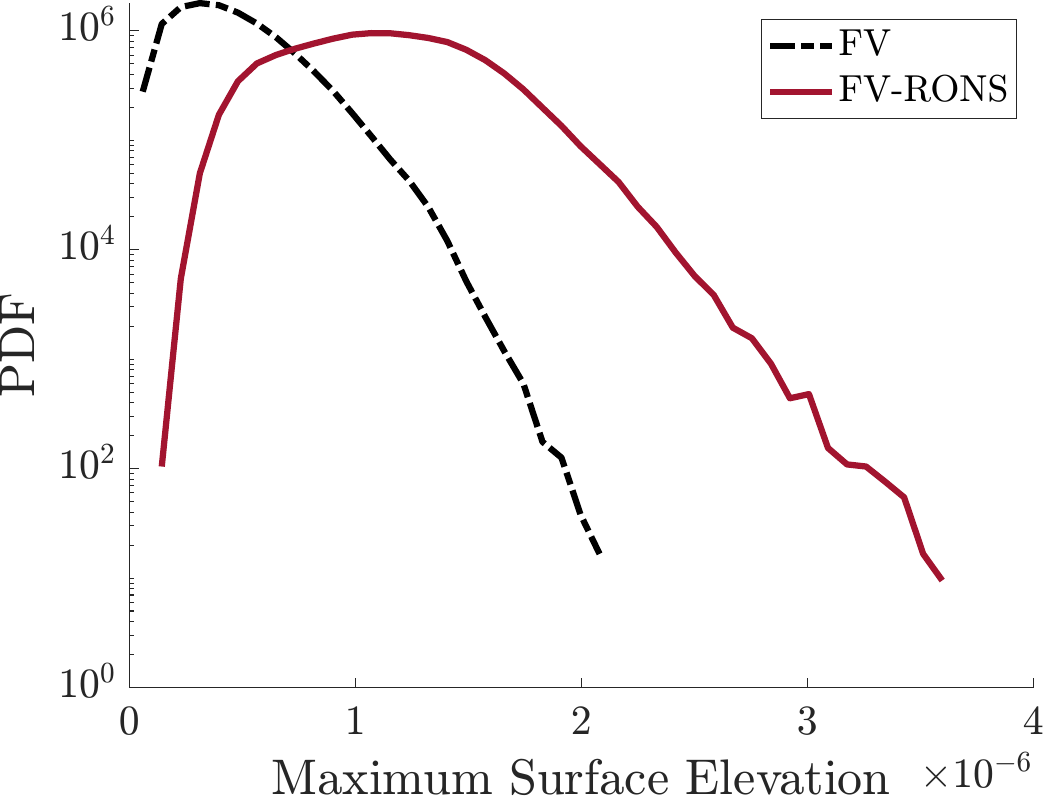}
\caption{\label{fig:fvrons_maxes} Comparisons of the PDFs of the maximum surface elevation $\max_{x}\eta(x,t)$ for traditional finite volume (dashed black line) and FV RONS (solid red line). The surface elevations are reported in dimensionless variables.}
\end{figure}

Figure \ref{fig:fvrons_maxes} compares the resulting PDFs for both traditional FV method and FV-RONS. As a result of energy decay, the FV method significantly underestimates the 
surface elevation. The most probable surface elevation according to the FV method is around $0.5\times 10^{-6}\simeq 1.1$~meter. In contrast, the FV-RONS has a mode at approximately $1.5\times 10^{-6}\simeq 3.2$~meters. Moreover, FV-RONS has a heavier tail as compared to the FV solution.
These results are in agreement with figure~\ref{fig:fvrons_etas}; as the energy dissipates, the wave heights decrease as well.

\subsection{Nonlinear Sch\"{o}dinger Equation}\label{sec:nls}
As a second example, we consider the nonlinear Schr\"{o}dinger equation, modeling unidirectional, slowly modulating surface waves in a deep fluid~\cite{Anderson2022b,cousins16,Farazmand2017}. 
Denoting the surface elevation by $\tilde\eta(\tilde x,\tilde t)$, we consider perturbations to a sinusoidal wave such that
\begin{equation}\tilde\eta(\tilde x,\tilde t) = \mbox{Re}\left\{\tilde u(\tilde x,\tilde t)\exp[\hat i(k_0\tilde x-\omega_0\tilde t)]\right\}.
\end{equation}
 where $\tilde u$ denotes the complex-valued wave envelope. Here, $k_0$ is the wave number of the sinusoidal carrier wave and $\omega_0$ is its angular frequency. Defining the dimensionless variables $t = \omega_0 \tilde{t}$, $x = k_0 \tilde{x}$, $u = k_0 \tilde{u}$, the complex envelope $u(x,t)$ satisfies the nonlinear Schr\"odinger equation~\cite{Zakharov68}, 
\begin{equation}\label{eq:NLS}
u_t =- \frac{1}{2}u_x -\frac{\hat{i}}{8}u_{xx}-\frac{\hat{i}}{2}|u|^2u.
\end{equation} 
We note that a similar equation has been used to describe dispersive optical waves~\cite{Onorato13}.

NLS admits a hierarchy of first integrals~\cite{tao2006}. In our reduced model, we enforce two of the physically most relevant conserved quantities,
\begin{subequations}\label{eq:NLS_cq}
\begin{equation}\label{eq:NLS_mass}
I_1 = \int_{\R}|u|^2 \id x , 
\end{equation}
\begin{equation}
I_2  = \frac{1}{8}\int_{\R}|u_x|^2 \id x - \frac{1}{4}\int_{\R}|u|^4 \id x,
\label{eq:NLS_energy}
\end{equation}
\end{subequations}
which represent the total mass and energy, respectively. 
 \begin{figure}
 	\centering
	\includegraphics[width=.75\textwidth]{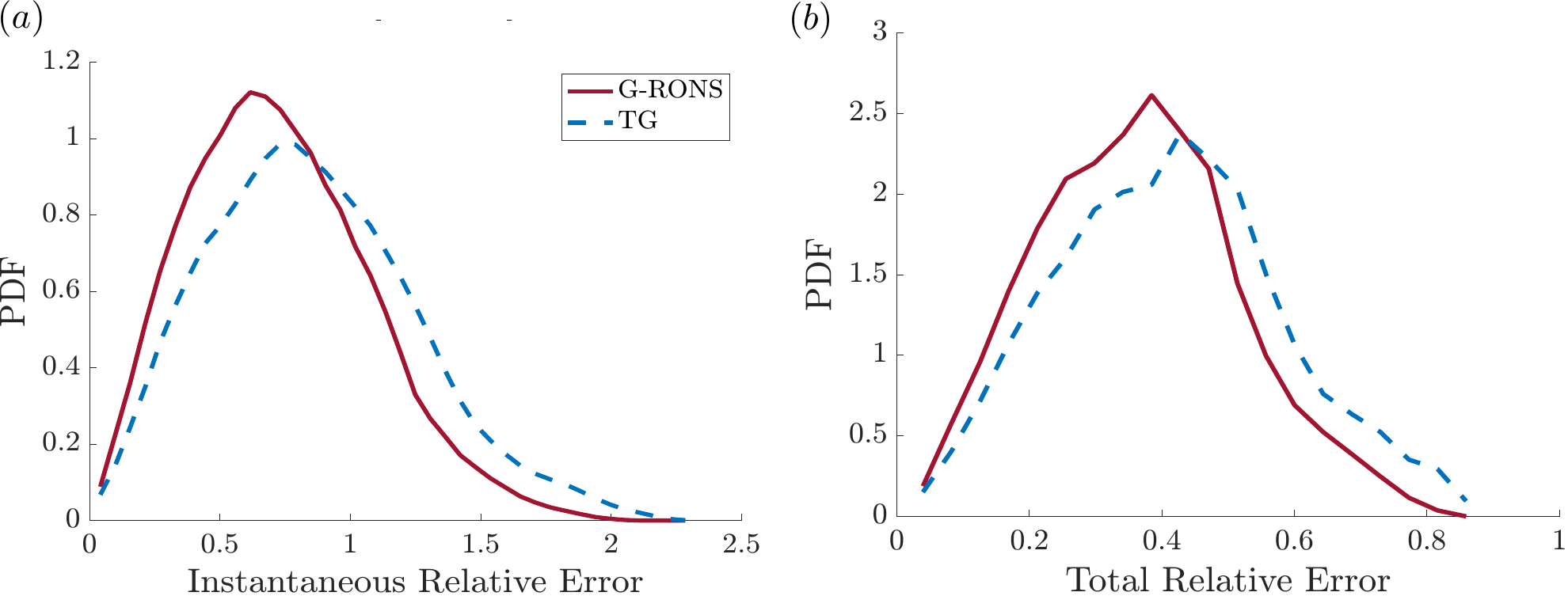}
	\caption{PDF of the relative errors: (a) Instantaneous relative error~\eqref{eq:inst_err}, (b) Total relative error~\eqref{eq:total_err} computed from $10^4$ initial conditions over the time interval $[400,600]$. In each panel the dashed blue line indicates the traditional Galerkin projection, while the solid red line denotes G-RONS.}
	\label{fig:NLS_errs}
\end{figure}

For the reduced-order model, we consider Galerkin-type projections onto POD modes. More specifically, we consider reduced-order solutions of the form,
\begin{equation}\label{eq:nls_pod}
\hat{u}(x,\vc a(t)) = \overline u(x)+ \sum_{i=1}^N a_i(t) \phi_i(x)
\end{equation}
where $\phi_i$ denotes the $i$-th POD mode, and $\overline u(x)$ is the mean. In order to compute the POD modes, we follows a standard algorithm as outlined in Ref.~\cite{Shlizerman2012}.
Here we consider NLS over a domain $x\in[0,L]$ with periodic boundary conditions. To gather data for the POD algorithm, we run multiple simulations of NLS over the time interval $t\in[0,2000\pi]$ (a thousand wave periods) using a standard pseudo-spectral method with $2^{10}$ Fourier modes for the spatial discretization. 
Each simulation corresponds to a different initial condition with random phases. More specifically, we define
\begin{equation}
\tilde{u}_0 (x)= \sum_{j=3}^8 e^{-\frac{j^2}{10}}\cos\left(\frac{2\pi  j x}{L} + \varphi_j\right),
\end{equation}
where the domain size $L = 256\pi$ is equivalent to 128 wavelengths, and $\varphi$ are the random phases drawn from a uniform distribution over $[0,2\pi]$. In order to ensure that the amplitude of the initial envelope $u_0$ is realistic, we rescale $\tilde{u}_0$ to define the initial condition,
\begin{equation}
u_0 (x)= \frac{0.13}{\max_x\{\tilde u_0(x)\}}\tilde{u}_0(x).
\end{equation}
 For the integration in time, we use an exponential time-differencing scheme as described in \cite{Cox2002} with a time step of $\Delta t = 0.025$ in accordance with \cite{Anderson2022b,cousins16}. 
 
We evolve the model parameters $a_i(t)$ using two methods: traditional Galerkin (TG) projection and Galerkin projection with RONS (G-RONS). In G-RONS, we ensure that the first integrals~\eqref{eq:NLS_cq} are conserved in the reduced model~\eqref{eq:rons_systemPDE}. In contrast, TG does not necessarily respect these conserved quantities. For both TG and G-RONS, we use $N=9$ POD modes and compare our results against DNS solutions. Note that, although NLS is a complex-valued PDE, it can be written as a system of two real-valued PDEs for the real and imaginary parts of the solution. Therefore, the method introduced in section~\ref{sec:RONS_systems} is applicable here.
\begin{figure}
	\centering 
	\includegraphics[width=.9\textwidth]{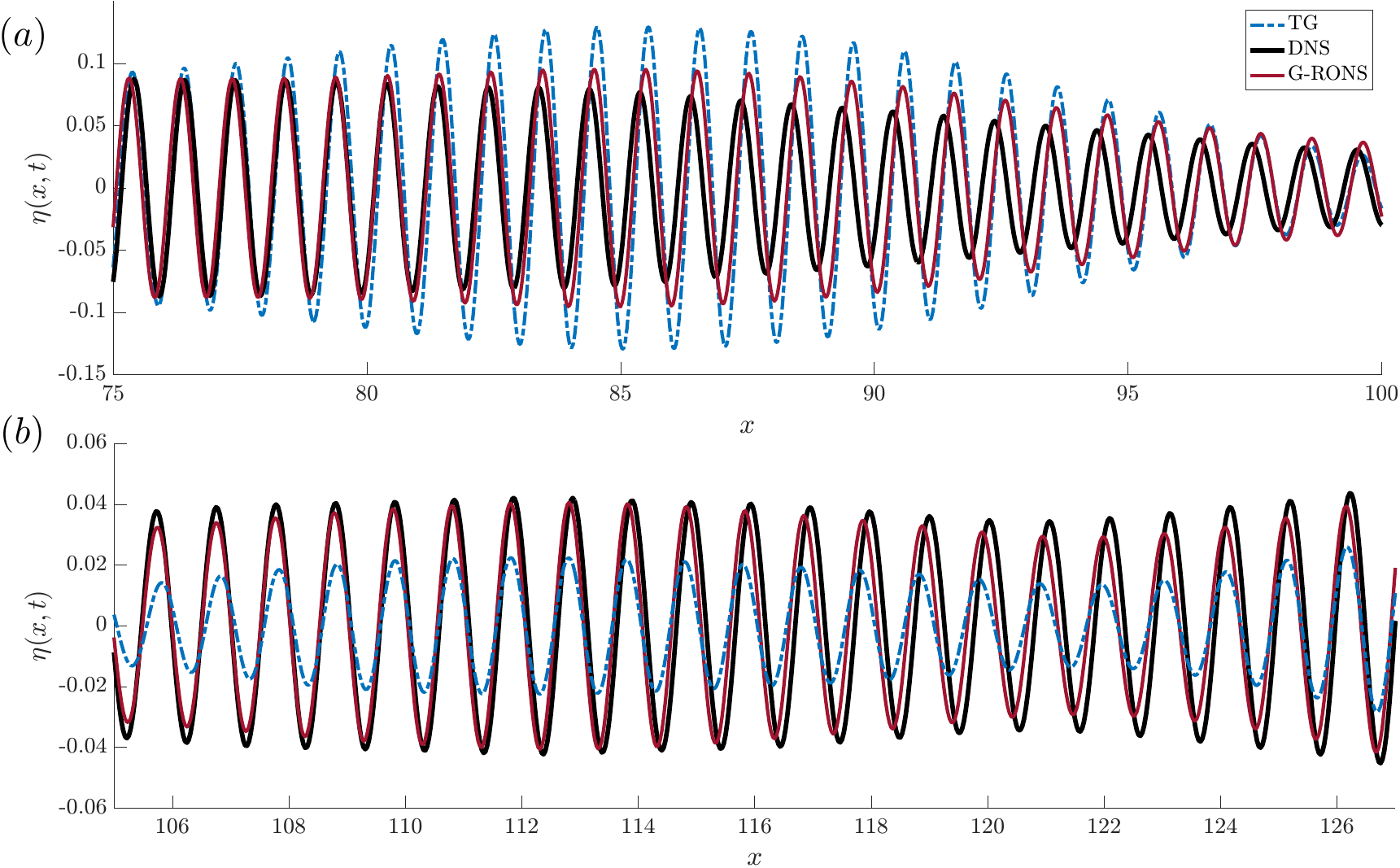}
	\caption{Two snapshot of the wave elevation $\eta(x,t) = \mbox{Re}[u(x,t)\exp(\hat{i}(x-t))]$ comparing the TG model (dashed blue), the G-RONS model (solid red) and the ground truth (solid black).}
	\label{fig:NLS_snap} 
\end{figure}

For the reduced models, we chose $10^4$ different initial conditions by sampling the initial amplitudes $a_i(0)$ of the first five POD modes ($1\leq i\leq 5$) from a standard uniform distribution over the interval $[0,1]$ and set the remaining initial amplitudes ($6\leq i\leq 9$) equal to zero. To compare the accuracy of the models, we define the instantaneous and total  relative errors, respectively,
\begin{align}
\varepsilon_{I}(t) &= \frac{\int_0^L |u(x,t) - \hat{u}(x,t)|^2\id x }{\int_0^L |u(x,t) |^2\id x } \label{eq:inst_err}, \\
\varepsilon_{T} &= \frac{\int_{t_i}^{t_f} \int_{0}^L |u(x,t) - \hat{u}(x,t)|^2\id x \id t}{\int_{t_i}^{t_f} \int_0^L |u(x,t) |^2\id x \id t} \label{eq:total_err}, 
\end{align}
where $u(x,t)$ is the DNS solution and $\hat{u}(x,t)$ is the reduced-order solution either from TG or from G-RONS. The solutions are integrated for $t_f = 600$ time units and the first $400$ time units are discarded to ensure initial transients do not affect the results; therefore, we have $t_i=400$.
The instantaneous errors are recorded every 1 time unit which leads to $10^4$ measurements of the total error and $2\times 10^6$ measurements of the instantaneous error.

Figure \ref{fig:NLS_errs} shows the PDFs of the instantaneous and total errors. In both cases, the mode of the error occurs at a lower value for G-RONS than TG. Moreover, we can see lower relative errors occur more frequently for G-RONS, while higher relative errors occur less frequently. Both results indicate that, by enforcing the conserved quantities~\eqref{eq:NLS_cq}, we gain additional accuracy in the reduced-order model. Figure \ref{fig:NLS_snap} shows the surfaces elevation at two time instances. In both instances, the true surface elevation is much closer to the G-RONS reduced-order solution, whereas TG either overestimates (top panel) or underestimates (bottom panel) the true wave height. 

\begin{figure}
\centering 
\includegraphics[width=0.49\textwidth]{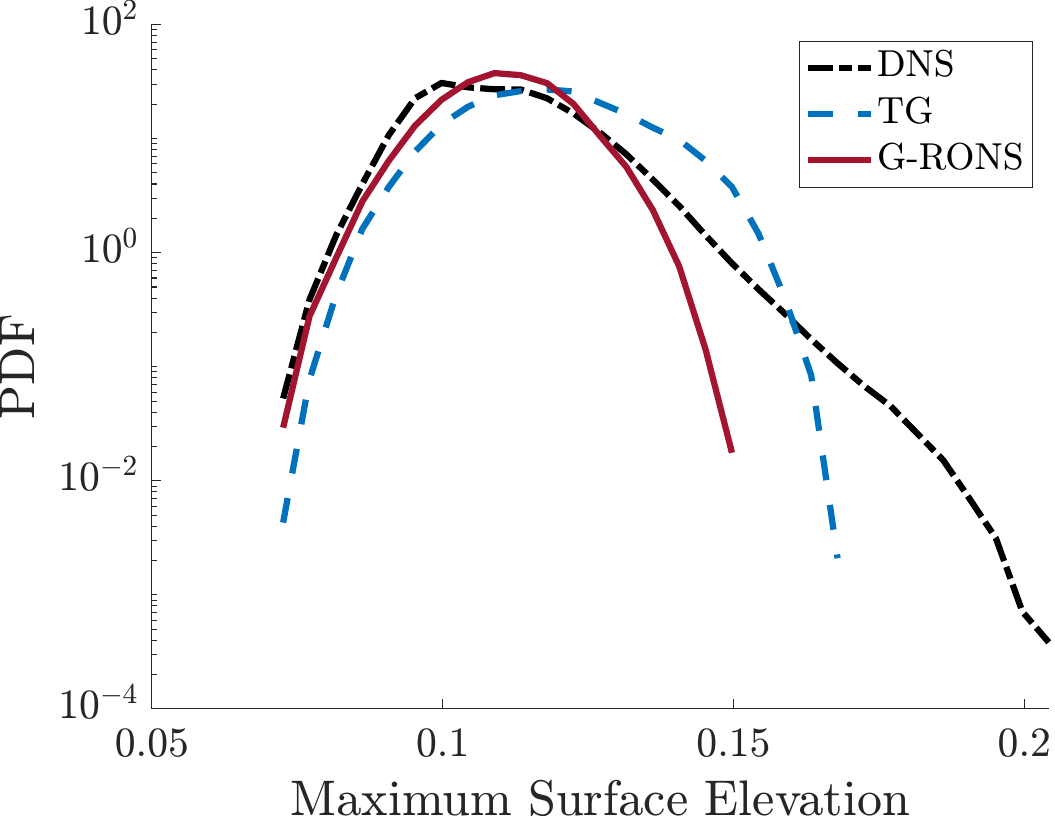}
\caption{\label{fig:NLS_max} PDF of the maximum envelope: $\max_{x}|u(x,t)|$. The PDFs are computed from data saved at every $\Delta t = 1$ time units in the interval $[200, 600]$ to discard the initial transients. DNS results are marked by the solid black line, TG is the dashed blue, and G-RONS is the solid red line.}
\end{figure}

Since NLS is frequently used in the prediction of rogue waves \cite{cousins16}, we also compute the PDF of the maximum surface elevation, i.e., $\max_{x} |u(x,t)|$. 
Figure \ref{fig:NLS_max} shows that both methods fail to correctly quantify the tail of the distribution corresponding to large surface elevations. However, the G-RONS model more closely reproduces the core of the PDF near its mode. In closing, we point out that Galerkin-type reduced-order models are known to underestimate the frequency of rogue waves (tail of the distribution). To address this issue, more complex nonlinear reduced-order models have been developed which depend nonlinearly on length-scale and phase parameters, in addition to the amplitude~\cite{cousins15,cousins16}. RONS has been shown to accurately evolve the parameters of such more complex models~\cite{Anderson2022a,Anderson2022b}. 

\section{Conclusions}\label{sec:conc}
It is well-known that certain properties of PDEs may be lost after their numerical discretization. For instance, first integrals of the PDE may not be invariant after
discretization. This issue has led to the development of a plethora of structure-preserving numerical schemes for direct numerical simulation and reduced-order modeling of PDEs (see section~\ref{sec:rel_work} for a review). However, these methods are often applicable to a particular PDE or a special class of PDEs such as Hamiltonian systems. 

Here, we developed two general-purpose methods, Galerkin RONS and finite volume RONS, which are applicable to evolution PDEs without any further assumption about the structure of the equations. Both methods are derived from appropriately formulated constrained optimization problems whose solutions can be explicitly obtained as a set of ODEs. The constraints ensure that the first integrals of the PDE remain conserved after discretization. 

The resulting ODEs are similar to those found from usual Galerkin or finite volume discretization with the crucial difference that a summation term, involving Lagrange multipliers, is added to the equations. As such, it is straightforward to implement our methods on top of existing Galerkin or finite volume code. 
At every time step, one needs to solve a linear algebraic system to obtain the Lagrange multiplier. The size of this system is small and equal to the number of conserved quantities being enforced. Therefore, the linear solves do not constitute a significant computational cost.

We demonstrated the application of our methods on two test cases. First, we considered a central-upwind finite volume schemes for the shallow water equation. Since this scheme does not preserve the total energy of the system, the resulting numerical solutions exhibit significant decay in the amplitude of the surface elevation. In contrast, finite volume RONS ensures the conservation of energy and consequently eliminates this spurious amplitude decay.

Our second example involves a POD-Galerkin reduced model of the nonlinear Schr\"odinger equation for deep water waves. Galerkin RONS ensures that the total mass and energy of the system are conserved. Comparing the results against direct numerical simulations, we observed that the Galerkin RONS reduced model is significantly more accurate than the usual POD-Galerkin model. 

Galerkin RONS~\eqref{eq:rons_eq} and finite volume RONS~\eqref{eq:rons_eq_fv} are guaranteed to conserve the first integrals $I_i$ in their continuous-time formulation. It is a possibility that this property may be lost after temporal discretization, e.g., due to numerical diffusion. In our numerical examples, we used standard Runge--Kutta schemes and did not observe any appreciable loss of conserved quantities. Nonetheless, conservative time integration schemes can be used to ensure that the conservative nature of the RONS equations survive their temporal discretization~\cite{gottlieb1998,Shu1988}.

As mentioned earlier, our methods are universally applicable to time-dependent PDEs. This universality comes at a cost: additional structures of the PDE, beyond its conserved quantities, may still be violated by Galerkin RONS and finite volume RONS. For instance, the two-form and the Poisson bracket associated with Hamiltonian systems may not survive our truncation. It would be interesting to investigate whether such additional structure can be incorporated in the RONS framework.

\section*{Acknowledgments}
We would like to thank Prof. Alina Chertock and Dr. Michael Reddle for fruitful conversations.
This work was supported by the National Science Foundation under the grants DMS-1745654 and DMS-2208541.

\appendix
\section{RONS with Complex Parameters}\label{app:rons_complex}
Even when the solution $u(\vc x,t)$ of the PDE is real-valued, the modes $\phi_i$ and the corresponding amplitudes $a_i$ might be complex-valued. This occurs, for instance, when using Fourier modes. In this appendix, we discuss how the RONS equations can be modified to allow for complex-valued modes and amplitudes. In this case, we decompose each amplitude $a_i$ into its real and imaginary parts; namely we let $a_i(t) = \alpha_i(t) + \hat{i}\beta_i(t)$. Then, the vector of parameters is defined as $\vc a = [\boldsymbol{\alpha}^\top \ \ \boldsymbol{\beta}^\top ]^\top \in \R^{2N}$, where $\bs \alpha$ and $\bs \beta$ are vectors whose $i$-th components are $\alpha_i$ and $\beta_i$, respectively. 

For the RONS optimization problem, the cost function remains unchanged, i.e. $\mathcal{J}(\vc a, \dot{\vc a}) = \frac{1}{2}\|\hat{u}_t - F(u)\|_H^2$. Thus the optimization problem reads
\begin{align}\label{eq:rons_const_opt_app}
\begin{split}
 \min_{\dot{\vc a}\in\R^{2N}} \  & \mathcal{J}(\vc a, \dot{\vc a}), \\
 \text{s.t.} \ & I_k(\vc a(t)) = \mbox{const.}, \quad \text{for }  k = 1,2,\ldots , m.
 \end{split}
\end{align}
For complex-valued parameters, the metric tensor $\vc M$ is constructed as follows. Let $\tilde{\vc M}$ have components $\tilde{M}_{ij} = \langle\phi_i,\phi_j \rangle_H$; then construct the metric tensor $\vc M$ as, 
\begin{equation}
\vc M = \begin{bmatrix}
\tilde{\vc M} &- \hat{i}\tilde{\vc M} \\ 
\hat{i}\tilde{\vc M} & \tilde{\vc M}
\end{bmatrix}.
\end{equation}
Moreover, let the vector field $\tilde{\vc f}(\vc a)$ has components $\tilde{f}_i = \langle \phi_i, F(\hat{u})\rangle_H$, and let $\vc f(\vc a) = [\tilde{\vc f}(\vc a)^\top \ \ \hat{i}\tilde{\vc f}(\vc a)^\top]^\top$. Then the minimizer of the unconstrained problem is 
\begin{equation}
\vc M^{(r)} \dot{\vc a} = \vc f^{(r)}(\vc a),
\end{equation}
where the superscript $(r)$ denotes the real part. The Lagrange multiplier $\boldsymbol{\lambda}$ is the solution to the constraint equation $\vc C \boldsymbol\lambda = \vc b$. The components of the constraint matrix and the vector $\vc b$ are defined similarly, with
\begin{align}\label{eq:constraint_complex}
\begin{split}
C_{ij} =\left \langle \nabla_{\vc a} I_i(\vc a),  \vc \left(\vc M^{(r)}\right)^{-1} \nabla_{\vc a} I_j (\vc a)\right \rangle, \\
b_i = \left \langle \nabla_{\vc a} I_i(\vc a) , \left(\vc M^{(r)}\right)^{-1} \vc f^{(r)}(\vc a)\right \rangle,
\end{split}
\end{align}
where the gradient $\nabla_{\vc a}I_j(\vc a) = [\nabla_{\bs\alpha} I_j(\vc a)^\top\ \nabla_{\vc \beta}I_j (\vc a)^\top]^\top$. The solution to the constrained minimization problem is 
\begin{equation}\label{eq:rons_complex}
\vc M^{(r)}\dot{\vc a} = \vc f^{(r)}(\vc a) - \sum_{j=1}^m \lambda_i \nabla_{\vc a} I_j(\vc a). 
\end{equation}

For a system of PDEs with $p$ components, we make the exact same adjustments to each set of $\vc a_i$, $\vc M_i$, and $\vc f_i$. Similar to section~\ref{sec:RONS_systems}, the metric tensor $\vc M$ is the block diagonal matrix, $$\vc M = \text{diag}(\vc M_1^{(r)}, \ldots, \vc M_p^{(r)}),$$ 
whose blocks are the real parts of 
\begin{equation}
\vc M_k = \begin{bmatrix}
\tilde{\vc M}_k & -\hat{i}\tilde{\vc M}_k\\
\hat{i}\tilde{\vc M}_k & \tilde{\vc M}_k\\ 
\end{bmatrix},
\end{equation}
and the vectors $\vc a$ and $\vc f$ are defined by
\begin{equation}
\vc a = \begin{bmatrix}
\bs \alpha_1 \\
\bs \beta_1 \\ 
\bs \alpha_2 \\ 
\bs \beta_2 \\
\vdots \\
\bs \alpha_p \\
\bs \beta_p \\
\end{bmatrix},
\quad \text{and}
\quad 
\vc f(\vc a) = \begin{bmatrix}
\tilde{\vc f_1}(\vc a)\\
\hat{i} \tilde{\vc f_1}(\vc a)\\
\tilde{\vc f_2}(\vc a)\\
\hat{i} \tilde{\vc f_2}(\vc a)\\
\vdots \\
\tilde{\vc f_p}(\vc a)\\
\hat{i} \tilde{\vc f_p}(\vc a)
\end{bmatrix}. 
\end{equation}
Furthermore, for systems of equations, the gradients of the conserved quantities are constructed similarly, taking 
\begin{equation}
\nabla_{\vc a} I_j(\vc a) = \begin{bmatrix}
\nabla_{\bs\alpha_1}I_j (\vc a)\\ 
\nabla_{\bs\beta_1}I_j(\vc a) \\ 
\nabla_{\bs\alpha_2}I_j (\vc a)\\ 
\nabla_{\bs\beta_2}I_j (\vc a)\\ 
\vdots \\
\nabla_{\bs\alpha_p}I_j(\vc a) \\ 
\nabla_{\bs\beta_p}I_j (\vc a)\\ 
\end{bmatrix}.
\end{equation}
Constructing the constraint matrix $\vc C$ and the vector $\vc b$ as in \eqref{eq:constraint_complex}, we arrive at the RONS equation in the same form as Eq. \eqref{eq:rons_complex}.

%\bibliographystyle{plain}      % mathematics and physical sciences
%\bibliography{biblio.bib,bibliog.bib}   % name your BibTeX data base

\end{document}